\documentclass[aap,MSNbibl,dvips]{arximspdf}

\doi{10.1214/12-AAP861} 
\volume{23}
\issue{3}
\pubyear{2013}
\firstpage{957}
\lastpage{986}

\makeatletter
\newtheorem{lemma}{Lemma}[section]
\newtheorem{theorem}[lemma]{Theorem}
\newtheorem{proposition}[lemma]{Proposition}

\newcommand{\implies}{\Longrightarrow}

\newcommand{\cal}{\mathcal}

\newcommand{\PP}{\mathbb{P}}
\newcommand{\ve}{{\varepsilon}}
\newcommand{\pd}{\partial}
\newcommand{\E}{{\mathbb E}}
\renewcommand{\O}{{\Omega}}
\newcommand{\R}{{\mathbb R}}
\newcommand{\N}{{\mathbb N}}
\newcommand{\Z}{{\mathbb Z}}
\renewcommand{\a}{{\alpha}}
\renewcommand{\b}{{\beta}}
\newcommand{\D}{{\Delta}}
\renewcommand{\t}{{\tau}}
\renewcommand{\k}{{\kappa}}
\newcommand{\g}{{\gamma}}
\newcommand{\s}{{\sigma}}
\renewcommand{\L}{{\Lambda}}
\newcommand{\cF}{{\cal F}}
\newcommand{\bX}{{\mathbf X}}

\newcommand{\sse}{\subseteq}

\makeatother

\begin{document}
\begin{frontmatter}

\title{Averaging over fast variables in the fluid limit for Markov chains:
Application to the supermarket model with memory}
\runtitle{Averaging over fast variables}

\begin{aug}
\author[A]{\fnms{M. J.} \snm{Luczak}\corref{}\thanksref{t1}\ead[label=e1]{m.luczak@qmul.ac.uk}}
\and
\author[B]{\fnms{J. R.} \snm{Norris}\thanksref{t2}\ead[label=e2]{j.r.norris@statslab.cam.ac.uk}}
\runauthor{M. J. Luczak and J. R. Norris}
\affiliation{London School of Economics and University of Cambridge}
\address[A]{School of Mathematical Sciences\\
Queen Mary, University of London\\
Mile End Road\\
London E1 4NS\\
United Kingdom\\
\printead{e1}} 
\address[B]{Statistical Laboratory\\
Centre for Mathematical Sciences\\
University of Cambridge\\
Wilberforce Road\\
Cambridge, CB3 0WB\\
United Kingdom\\
\printead{e2}}
\end{aug}

\thankstext{t1}{Supported in part by a STICERD grant at the LSE, and by EPSRC Leadership Fellowship EP/J004022/2.
Part of this work was accomplished while visiting the Mittag-Leffler
Institute.}

\thankstext{t2}{Supported by EPSRC Grant EP/E01772X/1.}

\received{\smonth{3} \syear{2010}}
\revised{\smonth{8} \syear{2011}}

%
\begin{abstract}
We set out a general procedure which allows the approximation of
certain Markov chains by the solutions of differential equations. The
chains considered have some components which oscillate rapidly and
randomly, while others are close to deterministic. The limiting
dynamics are obtained by averaging the drift of the latter with respect
to a local equilibrium distribution of the former. Some general
estimates are proved under a uniform mixing condition on the fast
variable which give explicit error probabilities for the fluid
approximation.
Mitzenmacher, Prabhakar and Shah [In \textit{Proc. 43rd Ann. Symp. Found.
Comp. Sci.} (2002) 799--808, IEEE] introduced a variant with memory of the
``join the shortest queue'' or ``supermarket'' model, and obtained a
limit picture for the case of a stable system in which the number of
queues and the total arrival rate are large. In this limit, the
empirical distribution of queue sizes satisfies a differential
equation, while the memory of the system oscillates rapidly and
randomly. We illustrate our general fluid limit estimate by giving a
proof of this limit picture.
\end{abstract}

%

\end{frontmatter}

\section{A general fluid limit estimate}\label{GFL}
We describe a general framework to allow the incorporation of averaging
over fast variables into
fluid limit estimates for Markov chains, building on the approach used
in~\cite{MR2395153}.
The main results of this section, Theorems~\ref{FLE} and~\ref{FLET},
establish explicit error probabilities
for the fluid approximation under assumptions which can be verified from
knowledge of the transition rates of the Markov chain. Also see \cite
{MR2288709} for related results.

\subsection{Outline of the method}\label{OM}
Let $X=(X_t)_{t\ge0}$ be a continuous-time Markov chain with countable
state-space $S$
and with generator matrix $Q=(q(\xi,\xi')\dvtx\break\xi,\xi'\in S)$.
Assume that the total jump rate $q(\xi)$ is finite for all states $\xi
$, and that $X$ is nonexplosive.
Then the law of $X$ is determined uniquely by $Q$ and the law of $X_0$.
Make a choice of \textit{fluid coordinates} $x^i\dvtx S\to\R$, for
$i=1,\ldots
,d$, and write ${\mathbf x}=(x^1,\ldots,x^d)\dvtx S\to\R^d$.
Consider the $\R^d$-valued process $\bX=(\bX_t)_{t\ge0}$ given by
$\bX_t=(X^1_t,\ldots,X^d_t)={\mathbf x}(X_t)$.
Call $\bX$ the \textit{slow} or \textit{fluid variable}.
Define for each $\xi\in S$ the \textit{drift vector}
\[
\b(\xi)=Q{\mathbf x}(\xi)=\sum_{\xi'\not=\xi}\bigl({\mathbf x}(\xi
')-{\mathbf x}(\xi)\bigr)q(\xi,\xi').
\]
Also, make a choice of an \textit{auxiliary coordinate} $y\dvtx S\to
I$, for
some countable set~$I$, and set $Y_t=y(X_t)$.
Call the process $Y=(Y_t)_{t\ge0}$ the \textit{fast variable}.
For $\xi\in S$ and $y'\in I$ with $y'\not=y(\xi)$, write $\g(\xi
,y')$ for the total rate at which $Y$ jumps to $y'$ when $X$ is at $\xi$.
Thus
\[
\g(\xi,y')=\sum_{\xi'\dvtx y(\xi')=y'}q(\xi,\xi').
\]
Choose a subset $U$ of $\R^d$ and a function $b\dvtx U\times I\to\R^d$.
Choose also, for each $x\in U$, a generator matrix
$G_x=(g(x,y,y')\dvtx y,y'\in I)$ having a unique
invariant distribution $\pi_x=(\pi(x,y)\dvtx y\in I)$.
These choices are to be made so that $\b(\xi)$ is close to
$b({\mathbf x}(\xi),y(\xi))$ and $\g(\xi,y')$ is close
to $g({\mathbf x}(\xi),y(\xi),y')$ whenever ${\mathbf x}(\xi)\in U$
and $y'\in I$.
Define for $x\in U$
\[
\bar b(x)=\sum_{y\in I}b(x,y)\pi(x,y).
\]
Then,\vspace*{1pt} under regularity assumptions to be specified later,
there exists a function $\chi\dvtx U\times I\to\R^d$ such that
%
%
\begin{equation}\label{GCH-0}
G\chi(x,y)=\sum_{y'\in I}g(x,y,y')\chi(x,y')=b(x,y)-\bar b(x).
\end{equation}
Make a choice of such a function $\chi$. Call $\chi$ the
\textit{corrector for $b$}.

Fix $x_0\in U$. We will assume that $\bar b$ is Lipschitz on $U$.
Then the differential equation $\dot x_t=\bar b(x_t)$ has a unique
maximal solution $(x_t)_{t<{\zeta}}$ in $U$ starting from $x_0$.
Fix $t_0\in[0,{\zeta})$.
Then for $t\le t_0$,
%
%
\begin{equation}\label{DE}
x_t=x_0+\int_0^t\bar b(x_s)\,ds.
\end{equation}
Define for $\xi\in S$ with ${\mathbf x}(\xi)\in U$
\[
\bar{\mathbf x}(\xi)={\mathbf x}(\xi)-\chi({\mathbf x}(\xi),y(\xi)).
\]
Let $T$ be a stopping time such that $\bX_t\in U$ for all $t\le T$.
Then, under regularity assumptions to be specified later, for $t\le T$,
%
%
\begin{equation}\label{ME}
\bar{\mathbf x}(X_t)=\bar{\mathbf x}(X_0)+M_t+\int_0^t\bar\b(X_s)\,ds,
\end{equation}
where $M=M^{\bar{\mathbf x}}$ is a martingale and where
%
%
\begin{equation}\label{BB}
\bar\b=Q\bar{\mathbf x}=\b-Q(\chi({\mathbf x},y)).
\end{equation}
On subtracting equations (\ref{DE}) and (\ref{ME}) we obtain for
$t\le T\wedge t_0$
%
%
\begin{eqnarray}\label{KEY}
\bX_t-x_t&=&\bX_0-x_0+\chi(\bX_t,Y_t)-\chi(\bX_0,Y_0)+M_t+\int
_0^t\D(X_s)\,ds\nonumber\\[-8pt]\\[-8pt]
&&{} +\int_0^t\bigl(\b(X_s)-b(\bX_s,Y_s)\bigr)\,ds+\int_0^t\bigl(\bar b(\bX
_s)-\bar b(x_s)\bigr)\,ds,\nonumber
\end{eqnarray}
where $\D=G\chi({\mathbf x},y)-Q(\chi({\mathbf x},y))$.

The discussion in the present paragraph is intended for orientation
only, and will play no essential role in the derivation of our results.
Fix $U_0\sse U$ such that for all $\xi,\xi'\in S$ with ${\mathbf
x}(\xi)\in U_0$ and $q(\xi,\xi')>0$ we have ${\mathbf x}(\xi')\in U$.
Assume that $T$ is chosen so that $\bX_t\in U_0$ for all $t\le T$.
Define for $\xi\in S$ with ${\mathbf x}(\xi)\in U_0$ the
\textit{diffusivity tensor} $\a(\xi)\in\R^d\otimes\R^d$ by
%
%
\begin{equation}\label{DTD}
\a^{ij}(\xi)=\sum_{\xi'\not=\xi}\bigl(\bar{\mathbf x}^i(\xi')-\bar
{\mathbf x}^i(\xi)\bigr)\bigl(\bar{\mathbf x}^j(\xi')-\bar{\mathbf x}^j(\xi
)\bigr)q(\xi,\xi')
\end{equation}
and define for $t\le T$
\[
N_t=M_t\otimes M_t-\int_0^t\a(X_s)\,ds.
\]
Then, under regularity assumptions, $N$ is a martingale in $\R
^d\otimes\R^d$.
Choose a function $a\dvtx U_0\times I\to\R^d\otimes\R^d$ and set
\[
\bar a(x)=\sum_{y\in I}a(x,y)\pi(x,y).
\]
This choice is to be made so that $\a(\xi)$ is close to $a({\mathbf
x}(\xi),y(\xi))$ whenever \mbox{${\mathbf x}(\xi)\in U_0$}.
Suppose we can also find a corrector for $a$, that is, a function
$\tilde\chi\dvtx U_0\times I\to\R^d\otimes\R^d$ such that
%
%
\begin{equation}\label{GCH-1}
G\tilde\chi(x,y)=a(x,y)-\bar a(x).
\end{equation}
Then, for $t\le T$,
%
%
\begin{eqnarray}\label{KEYT}
\int_0^t\a(X_s)\,ds&=&\tilde\chi(\bX_t,Y_t)-\tilde\chi(\bX
_0,Y_0)-\tilde M_t+\int_0^t\tilde\D(X_s)\,ds\nonumber\\[-8pt]\\[-8pt]
&&{} +\int_0^t\bigl(\a(X_s)-a(\bX_s,Y_s)\bigr)\,ds+\int_0^t\bar
a(\bX_s)\,ds,\nonumber
\end{eqnarray}
where $\tilde\D=G\tilde\chi({\mathbf x},y)-Q(\tilde\chi({\mathbf
x},y))$ and, under suitable regularity conditions, $\tilde M=M^{\tilde
\chi}$ is a martingale up to $T$.

The martingale terms $M$ and $\tilde M$ in (\ref{KEY}) and (\ref
{KEYT}) can be shown to be small, under
suitable conditions, using the following standard type of exponential
martingale inequality.
In the form given here it may be deduced, for example, from \cite
{MR2395153}, Proposition 8.8, by
setting $f=\theta\phi$, $A=\theta^2e^{\theta J}\ve/2$ and $B=\theta
{\delta}$.
%
%
\begin{proposition}\label{EMI}
Let $\phi$ be a function on $S$. Define
\[
M_t=M^\phi_t=\phi(X_t)-\phi(X_0)-\int_0^t Q\phi(X_s)\,ds.
\]
Write $J=J(\phi)$ for the maximum possible jump in $\phi(X)$, thus
\[
J=\sup_{\xi,\xi'\in S, q(\xi,\xi')>0}|\phi(\xi')-\phi(\xi)|.
\]
Define a function $\a=\a^\phi$ on $S$ by
\[
\a(\xi)=\sum_{\xi'\not=\xi}\{\phi(\xi')-\phi(\xi)\}^2q(\xi
,\xi').
\]
Then, for all ${\delta},\ve\in(0,\infty)$ and all stopping times
$T$, we have
\[
\PP\biggl(\sup_{t\le T}M_t\ge{\delta}\mbox{ and }\int_0^T\a
(X_t)\,dt\le\ve\biggr)\le\exp\{-{\delta}^2/(2\ve e^{\theta
J})\},
\]
where $\theta\in(0,\infty)$ is determined by $\theta e^{\theta
J}={\delta
}/\ve$.
\end{proposition}

Now, if $\b,\g,\a$ are well approximated by $b,g,a$ and if we can
show that the corrector terms in (\ref{KEY})
and (\ref{KEYT}) are insignificant, then we may hope to use these
equations to show that
the path $(x_t\dvtx t\le t_0)$ provides a good (first order) approximation
to $(\bX_t\dvtx t\le t_0)$ and, moreover, that
the fluctuation process $(\bX_t-x_t\dvtx t\le t_0)$ is approximated (to
second order) by a Gaussian process $(F_t\dvtx t\le t_0)$
given by
\[
F_t=F_0+B_t+\int_0^t\nabla\bar b(x_s)F_s\,ds,
\]
where $(B_t\dvtx t\le t_0)$ is a zero-mean Gaussian process in $\R^d$ with
covariance
\[
\E(B_s\otimes B_t)=\int_0^{s\wedge t}\bar a(x_r)\,dr.\vadjust{\goodbreak}
\]
Our aim in the rest of this section is to give an explicit form of the
first order approximation with optimal
error scale, that is, of the same order as the scale of deviation
predicted by the second order approximation.
The next subsection contains some preparatory material on correctors.
A reader who wishes to understand only the statement of the fluid limit
estimate can skip directly to Section~\ref{SOTE}.\vspace*{-2pt}

\subsection{Correctors}\label{EC}
In order to implement the method just outlined, it is necessary either
to come up with
explicit correctors or to appeal to a general result which guarantees
the existence,
subject to verifiable conditions, of correctors with good properties.
In this subsection we obtain such a general result. In fact, we shall
find conditions which
guarantee the existence, for each bounded measurable function $f$ on
$U\times I$,
of a good \textit{corrector for $f$}, that is to say, a function $\chi
=\chi_f$ on $U\times I$ such that
\[
G\chi(x,y)=f(x,y)-\bar f(x),
\]
where
\[
\bar f(x)=\sum_{y\in I}f(x,y)\pi(x,y).
\]
Moreover, we shall see that $\chi_f$ depends linearly on $f$ and we
shall obtain a uniform bound and a
continuity estimate for $\chi_f$.

Assume that there is a constant $\nu\in(0,\infty)$ such that, for
all $x\in U$ and all $y\in I$,
the total rate of jumping from $y$ under $G_x$ does not exceed $\nu$.
Then we can choose an auxiliary measurable space $E$, with a $\sigma
$-field ${\mathcal E}$,
a family of probability measures $\mu=(\mu_x\dvtx x\in U)$ on $(E,
{\mathcal E})$ and a measurable function
$F\dvtx I\times E\to I$ such that, for all $x\in U$ and all $y,y'\in I$
distinct,
%
%
\begin{equation}\label{GNM}
g(x,y,y')=\nu\mu_x\bigl(\{v\in E\dvtx F(y,v)=y'\}\bigr).
\end{equation}
Let $N=(N(t)\dvtx t\ge0)$ be a Poisson process of rate $\nu$.
Fix $x\in U$ and let $V=(V_n\dvtx n\in\N)$ be a sequence of independent
random variables in $E$, all with law $\mu_x$.
Thus
\[
g(x,y,y')=\nu\PP\bigl(F(y,V_n)=y'\bigr)
\]
for all pairs of distinct states $y,y'$ and all $n$.
Fix a reference state $\bar y\in I$.
Given $y\in I$, set $Z_0=y$ and $\bar Z_0=\bar y$ and define
recursively for $n\ge0$,
\[
Z_{n+1}=F(Z_n,V_{n+1}),\qquad \bar Z_{n+1}=F(\bar Z_n,V_{n+1}).
\]
Set $Y_t=Z_{N(t)}$ and $\bar Y_t=\bar Z_{N(t)}$. Then $Y=(Y_t)_{t\ge
0}$ and $\bar Y=(\bar Y_t)_{t\ge0}$
are both Markov chains in $I$ with generator matrix $G_x$, starting
from $y$ and $\bar y$, respectively,\setcounter{footnote}{2}\footnote{The
process $Y$ introduced here is not the fast variable, also denoted $Y$
in the rest of the paper: the current $Y$ is
to be considered as a local approximation of the fast variable.\label{NFV}}
and are realized on the same probability space. We call the triple
$(\nu,\mu,F)$ a \textit{coupling mechanism}.
Define the \textit{coupling time}
\[
T_c=\inf\{t\ge0\dvtx Y_t=\bar Y_t\}.\vadjust{\goodbreak}
\]

Assume that, for some positive constant $\t$, for all $x\in U$ and all
$y,\bar y\in I$,
%
%
\begin{equation}\label{MTC}
m(x,y,\bar y)=\E_{(x,y,\bar y)}(T_c)\le\t.
\end{equation}
Fix a bounded measurable function $f$ on $U\times I$ and set
%
%
\begin{equation}\label{CFF}
\chi(x,y)=\E_{(x,y)}\int_0^{T_c}\bigl(f(x,Y_t)-f(x,\bar Y_t)\bigr)\,dt.
\end{equation}
Then $\chi$ is well defined and, for all $x\in U$ and all $y\in I$,
%
%
\begin{equation}\label{CTB}
|\chi(x,y)|\le2\t\|f\|_\infty.
\end{equation}

%
\begin{proposition}
The function $\chi$ is a corrector for $f$.
\end{proposition}
\begin{pf}
In the proof we suppress the variable $x$.
Note first that, if instead of taking $\bar Z_0=\bar y$, we start $\bar
Z$ randomly with the invariant distribution $\pi$, then
we change the value of $\chi$ by a constant independent of $y$. Hence,
it will suffice to establish the corrector equation
$G\chi=f-\bar f$ in this case.
Fix ${\lambda}>0$ and define
\[
\phi^{\lambda}(y)=\E\int_0^{T_{\lambda}}f(Y_t)\,dt,\qquad
\bar\phi^{\lambda}=\E\int_0^{T_{\lambda}}f(\bar Y_t)\,dt,
\]
where $T_{\lambda}=T_1/{\lambda}$, with $T_1$ an independent
exponential random variable of parameter $1$.
Then, since $Y$ and $\bar Y$ coincide after $T_c$,
\[
\bar\phi^{\lambda}-\phi^{\lambda}(y)
=\E\int_0^{T_{\lambda}\wedge T_c}\bigl(f(\bar Y_t)-f(Y_t)\bigr)\,dt\to\chi(y)
\]
as ${\lambda}\to0$.
By elementary conditioning arguments, $(G-{\lambda})\phi^{\lambda
}+f=0$ and \mbox{${\lambda}\bar\phi^{\lambda}=\bar f$}, so
\[
(G-{\lambda})(\bar\phi^{\lambda}-\phi^{\lambda})=f-\bar f.
\]
On passing to the limit ${\lambda}\to0$ in this equation, using
bounded convergence we find that $G\chi=f-\bar f$, as required.
\end{pf}

We remark that the corrector $\chi(x,\cdot)$ in fact depends only on $f$,
$G_x$ and the choice of $\bar y$,
as the preceding proof makes clear.
The further choice of a coupling mechanism is a way to obtain estimates
on $\chi$.

The following estimate will be used in dealing with the $\D$ term in
(\ref{KEY}).
We write $\|\mu_x-\mu_{x'}\|$ for the total variation distance
between $\mu_x$ and $\mu_{x'}$.
%
%
\begin{proposition}\label{CE}
For all $x,x'\in U$ and all $y\in I$,
%
%
\begin{eqnarray}\label{CEE}
|\chi(x,y)-\chi(x',y)|
&\le&
2\t\sup_{z\in I}|f(x,z)-f(x',z)|\nonumber\\[-8pt]\\[-8pt]
&&{}+2\nu\t
^2\|f\|_\infty\|\mu_x-\mu_{x'}\|.\nonumber
\end{eqnarray}
\end{proposition}
\begin{pf}
By a standard construction (maximal coupling) there exists a sequence
of independent random variables $((V_n,V_n')\dvtx n\in\N)$ in $E\times
E$ such that $V_n$ has distribution $\mu_x$, $V_n'$ has distribution
$\mu _{x'}$ and $\PP(V_n\not=V_n')=\frac12\|\mu_x-\mu_{x'}\|= \sup_{A
\in\mathcal E}|\mu_x(A)- \mu_{x'}(A)|$, for all $n$. Write
$(\cF_t)_{t\ge0}$ for the filtration of the marked Poisson process
obtained by marking $N$ with the random variables
$(V_n,V_n')$.\vspace*{2pt} Construct $(Y,\bar Y)$ from $N$ and
$(V_n\dvtx n\in\N)$ as above. Similarly\vspace*{1pt} construct
$(Y',\bar Y')$ from $N$ and $(V_n'\dvtx n\in\N)$. Recall that
$T_c=\inf\{t\ge0\dvtx Y_t=\bar Y_t\}$ and set $T_c'=\inf\{ t\ge0\dvtx
Y_t'=\bar Y_t'\}$. Set ${\lambda}=\frac12\nu\|\mu_x-\mu_{x'}\|$ and set
\[
D=\inf\{t\ge0\dvtx(Y_t,\bar Y_t)\not=(Y_t',\bar Y_t')\}.
\]
Then the process $t\mapsto1_{\{D\le t\}}-{\lambda}t$ is an $(\cF
_t)_{t\ge0}$-supermartingale
and $T_c$ is an $(\cF_t)_{t\ge0}$-stopping time.
So, by optional stopping, we have $\PP(D\le T_c)\le{\lambda}\E
(T_c)\le{\lambda}\t$.
Moreover, by the strong Markov property, on $\{D\le T_c\}$, we have $\E
(T_c-D|\cF_D)=m(x,Y_D,\bar Y_D)\le\t$ so,
for any function $g\dvtx I\to\R^d$, with $|g|\le\|f\|_\infty$,
\[
\E\biggl|\int_{D\wedge T_c}^{T_c}g(Y_t)\,dt\biggr|\le\t\|f\|_\infty
\PP(D\le T_c)\le{\lambda}\t^2\|f\|_\infty.
\]
On the other hand,
\[
\int_0^{D\wedge T_c}g(Y_t)\,dt=\int_0^{D\wedge T_c'}g(Y_t')\,dt
\]
so
\[
\biggl|\E\int_0^{T_c}g(Y_t)\,dt-\E\int_0^{T_c'}g(Y_t')\,dt\biggr|\le
2{\lambda}\t^2\|f\|_\infty=\nu\t^2\|f\|_\infty\|\mu_x-\mu_{x'}\|.
\]
We apply this estimate with $g=f(x,\cdot)$ to obtain
\begin{eqnarray*}
&&|\chi(x,y)-\chi(x',y)|\\
&&\qquad=\biggl|\E\int_0^{T_c}\bigl(f(x,\bar Y_t)-f(x,Y_t)\bigr)\,dt-\E\int
_0^{T_c'}\bigl(f(x',\bar Y_t')-f(x',Y_t')\bigr)\,dt\biggr|\\
&&\qquad\le 2\t\sup_{z\in I}|f(x,z)-f(x',z)|+\biggl|\E\int_0^{T_c}f(x,\bar
Y_t)\,dt-\E\int_0^{T_c'}f(x,\bar Y_t')\,dt\biggr|\\
&&\qquad\quad{} +\biggl|\E\int_0^{T_c}f(x,Y_t)\,dt-\E
\int_0^{T_c'}f(x,Y_t')\,dt\biggr|\\
&&\qquad\le 2\t\sup_{z\in I}|f(x,z)-f(x',z)|+2\nu\t^2\|f\|_\infty\|\mu
_x-\mu_{x'}\|
\end{eqnarray*}
as required.
\end{pf}

To summarize, we have shown the following proposition.
%
%
\begin{proposition}\label{GCR}
Assume conditions (\ref{GNM}) and (\ref{MTC}). Then, for any
bound\-ed measurable function $f$ on $U\times I$, there exists a corrector
$\chi_f$ for $f$ satisfying the estimates (\ref{CTB}) and (\ref{CEE}).
\end{proposition}

\subsection{Statement of the estimates}\label{SOTE}
Recall the context of Section~\ref{OM}. We consider a
continuous-time Markov chain $X$ with
countable state-space $S$ and generator matrix $Q$. We choose fluid
coordinates ${\mathbf x}\dvtx S\to\R^d$ and an
auxiliary coordinate $y\dvtx S\to I$. We choose also a subset $U\sse\R^d$,
which provides a means of
localization, together with a map $b\dvtx U\times I\to\R^d$, and a family
$G=(G_x\dvtx x\in U)$ of
generator matrices on $I$, each having a unique invariant distribution
$\pi_x$. Also choose, as
in the preceding subsection,
a coupling mechanism for $G$. This comprises a constant $\nu>0$, an
auxiliary space $E$,
a function $F\dvtx I\times E\to I$ and a family of probability distributions
$\mu=(\mu_x\dvtx x\in U)$ on $E$
such that
\[
g(x,y,y')=\nu\mu_x\bigl(\{v\in E\dvtx F(y,v)=y'\}\bigr),
\qquad x\in U, y,y'\in I\mbox{ distinct}.
\]
Define for $x\in U$
\[
\bar b(x)=\sum_{y\in I}b(x,y)\pi(x,y).
\]
Write $\bX_t={\mathbf x}(X_t)$ and assume that $(x_t)_{0\le t\le t_0}$
is a solution in $U$ to $\dot x_t=\bar b(x_t)$.
We use a scaled supremum norm on $\R^d$: fix positive constants $\s
_1,\ldots,\s_d$ and define for $x\in\R^d$
\[
\|x\|=\max_{1\le i\le d}|x_i|/\s_i.
\]

We now introduce some constants $\L,B,\t,J,J_1(b),J(\mu),K$ which
characterize
certain regularity properties of $Q$, $b$ and $G$.
Assume that, for all $\xi\in S$, all $x\in U$ and all $y,y'\in I$,
%
%
\begin{equation}\label{LBT}
q(\xi)\le\L,\qquad \|b(x,y)\|\le B,\qquad m(x,y,y')\le\t.
\end{equation}
Here $m(x,y,y')$ is the mean coupling time for $G_x$ starting from $y$
and $y'$, defined
in the preceding subsection, which depends on the choice of coupling mechanism.
Write $\cal J$ for the set of pairs of points in $U$ between which $\bX
$ can jump, thus
\[
{\cal J}=\{(x,x')\in U\times U\dvtx x={\mathbf x}(\xi),x'={\mathbf
x}(\xi
')\mbox{ for some $\xi,\xi'\in S$ with $q(\xi,\xi')>0$}\}.
\]
Set
\begin{eqnarray*}
J&=&\sup_{(x,x')\in{\cal J}}\|x-x'\|,\\
J_1(b)&=&\sup_{(x,x')\in{\cal J}, y\in I}\|b(x,y)-b(x',y)\|,\\
J(\mu)&=&\sup_{(x,x')\in{\cal J}}\|\mu_x-\mu_{x'}\|.
\end{eqnarray*}
Write $K$ for the Lipschitz constant of $\bar b$ on $U$; thus, for all
$x,x'\in U$,
%
%
\begin{equation}\label{BARB}
\|\bar b(x)-\bar b(x')\|\le K\|x-x'\|.
\end{equation}
Recall from Section~\ref{OM} the definitions of the drift vector
$\b$ for ${\mathbf x}$ and the jump rate $\g$ for $y$.
Define
\[
T=\inf\{t\ge0\dvtx\bX_t\notin U\}.
\]
Fix constants ${\delta}(\b,b),{\delta}(\g,g)\in(0,\infty)$ and
consider the events
%
%
\begin{equation}\label{OBB}
\O(\b,b)=\biggl\{\int_0^{T\wedge t_0}\|\b(X_t)-b({\mathbf
x}(X_t),y(X_t))\|\,dt\le{\delta}(\b,b)\biggr\}
\end{equation}
and
%
%
\begin{equation}\label{OGG}
\O(\g,g)=\biggl\{\int_0^{T\wedge t_0}\sum_{y'\not=y(X_t)}|\g
(X_t,y')-g({\mathbf
x}(X_t),y(X_t),y')|\,dt\le{\delta}(\g,g)\biggr\}.\hspace*{-35pt}
\end{equation}

%
\begin{theorem}\label{FLE}
Let $\ve>0$ be given and set ${\delta}=\ve e^{-Kt_0}/7$.
Assume that $J\le\ve$ and
\[
\max\bigl\{\|\bX_0-x_0\|, {\delta}(\b,b), 2\t B{\delta}(\g ,g), 2\t B,
2\L t_0\bigl(\t J_1(b)+\nu\t^2BJ(\mu)\bigr)\bigr\}\le{\delta}.
\]
Set $\bar J=J+4\t B$ and assume that ${\delta}\le\L\bar Jt_0/4$.
Further assume that the following \textit{tube condition} holds:
\[
\mbox{for }\xi\in S\mbox{ and }t\le t_0\qquad \|{\mathbf x}(\xi)-x_t\|\le
2\ve\quad\implies\quad{\mathbf x}(\xi)\in U.
\]
Then
\[
\PP\Bigl(\sup_{t\le t_0}\|\bX_t-x_t\|>\ve\Bigr)\le
2de^{-{\delta}^2/(4\L\bar J^2t_0)}+\PP\bigl(\O(\b,b)^c\cup\O(\g,g)^c\bigr).
\]
\end{theorem}

The proof of this result follows the initial stages of the proof of the
more elaborate Theorem~\ref{FLET} below.
We will not write it out separately but give further indications
immediately before the statement of Theorem~\ref{FLET}.
The reader will understand clearly the role of the inequalities which
appear as hypotheses
by following the proof. Here is an informal guide to their meanings.
The tube condition, together with $J\le\ve$,
allows us to localize the other hypotheses to $U$ by trapping the
process inside a tube around the limit path;
these conditions can be satisfied by choosing $U$ sufficiently large.
The conditions $\|\bX_0-x_0\|\le{\delta}$ and ${\delta}(\b,b)\le
{\delta}$
enforce that the initial conditions and drift fields match closely.
This requires, in particular, that
the fluid and auxiliary coordinates provide sufficient information to
nearly determine $\b$. The condition on ${\delta}(\g,g)$
forces a close match between the local behavior of the fast variable
and the idealized fast process used to compute the corrector.
The condition $2\t B\le{\delta}$ allows us to control the size of the
corrector, balancing the mean recurrence time
of the fast variable $\t$ against the range of the drift field $b$.
The condition on $2\L t_0(\t J_1(b)+\nu\t^2BJ(\mu))$
is needed for local regularity of the corrector, allowing us to pass
back from the idealized fast process at one
point $x$ to the actual fast variable when the fluid variable is near
$x$. Finally, the condition ${\delta}\le\L\bar Jt_0/4$
ensures we are in the ``Gaussian regime'' of the exponential martingale
inequality, where bad events cannot occur
by a small number of large jumps. For a nontrivial limiting dynamics,
$\L J$ should be of order $1$, while for
a useful estimate $\L J^2$ should be small; thus, as expected, we can
attempt to use the result when the Markov
chain takes small jumps at a high rate.

It is sometimes possible to improve on the constant $\L\bar J^2$
appearing in the preceding estimate, thereby
obtaining useful probability bounds for smaller choices of $\ve$.
However, to do this we have to make hypotheses expressed in terms of a
corrector.
Fix $\bar y\in I$ and denote by $\chi$ the corrector for $b$ given by
(\ref{CFF}).
Define for $\xi\in S$ with ${\mathbf x}(\xi)\in U$
\[
\bar{\mathbf x}(\xi)={\mathbf x}(\xi)-\chi({\mathbf x}(\xi),y(\xi)).
\]
Define, for $\xi\in S$ such that ${\mathbf x}(\xi)\in U$ and
${\mathbf x}(\xi')\in U$ whenever $q(\xi,\xi')>0$,
\[
\a^i(\xi)=\sum_{\xi'\not=\xi}\{\bar{\mathbf x}^i(\xi')-\bar
{\mathbf x}^i(\xi)\}^2q(\xi,\xi'),\qquad i=1,\ldots,d.
\]
Note that, since we shall be interested only in upper bounds, we deal
here only with the
diagonal terms of the diffusivity tensor defined at (\ref{DTD}).
Choose functions $a^i\dvtx I\to[0,\infty)$ such that, for all $\xi\in U$
where $\a^i(\xi)$ is defined,
%
%
\begin{equation}\label{ALA}
\a^i(\xi)\le a^i(y(\xi)),\qquad i=1,\ldots,d.
\end{equation}
For simplicity, we do not allow $a$ to depend on the fluid variable
${\mathbf x}(\xi)$. Since we can localize our
hypotheses near the (compact) limit path, we do not expect to lose much
precision by this simplification.
On the other hand, by permitting a dependence on the fast variable we
can sometimes do significantly
better than Theorem~\ref{FLE}, as we shall see in Section~\ref{SMM}. Set
\[
\bar a(x)=\sum_{y\in I}a(y)\pi(x,y),\qquad x\in U.
\]
We introduce two further constants $A$ and $\bar A$, with $\bar A\le
A\le\L\bar J^2$.
Assume that, for all $x\in U$ and all $y\in I$,
%
%
\begin{equation}\label{LBTT}
a^i(y)\le A\s^2_i,\qquad \bar a^i(x)\le\bar A\s_i^2,\qquad i=1,\ldots,d.
\end{equation}
Note\vspace*{1pt} that the corrector bound (\ref{CTB}) gives $\|\chi({\mathbf
x}(\xi),y(\xi))\|\le2\t B$, so $\a^i(\xi)\le\L\bar J^2\s_i^2$
and so (\ref{LBTT}) holds with $A = \bar A = \L\bar J^2$ and $a^i (y)
= A \sigma_i^2$ and $\bar a^i(x) = A \sigma_i^2$.
Thus Theorem~\ref{FLE} follows directly from (\ref{FLEA}) below.
The new inequalities required on the left-hand side of (\ref{ABAB}) can be
understood roughly as imposing that
the ratio of the averaged diffusivity to a uniform bound on the
diffusivity is not too small compared to
the mean recurrence time of the fast variable; so an effective
averaging takes place.
%
%
\begin{theorem}\label{FLET}
Assume that the hypotheses of Theorem~\ref{FLE} hold and that
${\delta}\bar J\le At_0/4$. Then
%
%
\begin{equation}\label{FLEA}
\PP\Bigl(\sup_{t\le t_0}\|\bX_t-x_t\|>\ve\Bigr)\le
2de^{-{\delta}^2/(4At_0)}+\PP\bigl(\O(\b,b)^c\cup\O(\g,g)^c\bigr).
\end{equation}
Moreover, under the further conditions ${\delta}\bar J\le\bar At_0/4$ and
%
%
\begin{equation}\label{ABAB}
\frac{1}{t_0}\max\{\t, \t{\delta}(\g,g), \L t_0\nu\t
^2J(\mu)\}\le\bar A/(20A)\le\L\t,
\end{equation}
we have
%
%
\begin{eqnarray}\label{FLEAB}
\PP\Bigl(\sup_{t\le t_0}\|\bX_t-x_t\|>\ve\Bigr)&\le&
2de^{-{\delta}^2/(4\bar At_0)}+2de^{-(\bar A/A)^2t_0/(6400\L\t
^2)}\nonumber\\[-8pt]\\[-8pt]
&&{}+\PP\bigl(\O(\b,b)^c\cup\O(\g,g)^c\bigr).\nonumber
\end{eqnarray}
\end{theorem}
\begin{pf}
Consider the stopping time
\[
T_0=\inf\{t\ge0\dvtx\|\bX_t-x_t\|>\ve\}.
\]
By the tube condition, we have $T_0\le T$. Moreover, for any $t<T_0$
and any $\xi'\in S$ such that $q(X_t,\xi')>0$,
we have
\[
\|{\mathbf x}(\xi')-x_t\|\le J+\|\bX_t-x_t\|\le2\ve
\]
so by the tube condition ${\mathbf x}(\xi')\in U$.

Recall\vspace*{1pt} that $\chi$ is the corrector for $b$ given by (\ref{CFF}).
For the proof of (\ref{FLEAB}), we shall use (\ref{CFF}) to also
construct a corrector $\tilde\chi$ for $a$.
Set $\tilde{\delta}=\bar At_0/10$.
Note from (\ref{CTB}) the bounds
\[
\|\chi(x,y)\|\le2\t B\le{\delta},\qquad |\tilde\chi^i(x,y)|\le2\t
A\s_i^2\le\tilde{\delta}\s_i^2.
\]
The inequality involving $\tilde{\delta}$ and further such
inequalities below, which depend on the first inequality in
assumption (\ref{ABAB}), will not be used in the proof of (\ref{FLEA}).
Write $\D=G\chi({\mathbf x},y)-Q(\chi({\mathbf x},y))=\D_1+\D_2$
and $\tilde\D=G\tilde\chi({\mathbf x},y)-Q(\tilde\chi({\mathbf
x},y))=\tilde\D_1+\tilde\D_2$,
where
%
%
\begin{equation}\label{DGH}
\D_1(\xi)=\sum_{y'\not=y(\xi)}\{g({\mathbf x}(\xi),y(\xi),y')-\g
(\xi,y')\}\chi({\mathbf x}(\xi),y')
\end{equation}
and
%
%
\begin{equation}\label{DGH-1}
\D_2(\xi)=\sum_{\xi'\not=\xi}q(\xi,\xi')\{\chi({\mathbf x}(\xi
),y(\xi'))-\chi({\mathbf x}(\xi'),y(\xi'))\}
\end{equation}
and where $\tilde\D_1$ and $\tilde\D_2$ are defined analogously.
Then, on $\O(\g,g)$, for $t\le T\wedge t_0$,
\[
\biggl\|\int_0^t\D_1(X_s)\,ds\biggr\|\le2\t B{\delta}(\g,g)\le
{\delta}
\]
and, using Proposition~\ref{CE},
\[
\biggl\|\int_0^t\D_2(X_s)\,ds\biggr\|\le2\L t_0\bigl(\t J_1(b)+\nu\t
^2BJ(\mu)\bigr)\le{\delta}.
\]
Similarly, for $t\le T\wedge t_0$,
\[
\biggl|\int_0^t\tilde\D^i_1(X_s)\,ds\biggr|\le2\t A{\delta}(\g
,g)\s_i^2\le\tilde{\delta}\s_i^2
\]
and
\[
\biggl|\int_0^t\tilde\D_2^i(X_s)\,ds\biggr|\le2\L t_0\nu\t^2AJ(\mu
)\s_i^2\le\tilde{\delta}\s_i^2.
\]
Take $M=M^{\bar{\mathbf x}}$ as in equations (\ref{ME}) and (\ref
{KEY}) and consider the event
\[
\O(M)=\Bigl\{\sup_{t\le T_0\wedge t_0}\|M_t\|\le{\delta}\Bigr\}.
\]
Then, on $\O(\b,b)\cap\O(\g,g)\cap\O(M)$, we can estimate the
terms in (\ref{KEY})
to obtain for $t\le T_0\wedge t_0$,
\[
\|\bX_t-x_t\|\le7{\delta}+K\int_0^t\|\bX_s-x_s\|\,ds,
\]
so that $\|\bX_t-x_t\|\le\ve$ by Gronwall's lemma. Note that this
forces $T_0\ge t_0$ and hence,
$\sup_{t\le t_0}\|\bX_t-x_t\|\le\ve$.
Set $\rho=3\bar A/2$ and consider the event
\[
\O(a)=\biggl\{\int_0^{T_0\wedge t_0}a^i(Y_s)\,ds\le\rho t_0\s_i^2
\mbox{ for } i=1,\ldots,d\biggr\}.
\]
By condition (\ref{ALA}), on $\O(a)$ we have
\[
\int_0^{T_0\wedge t_0}\a^i(X_s)\,ds\le\rho t_0\s_i^2.
\]
Set
\[
J_i=J(\bar{\mathbf x}^i)=\sup_{\xi,\xi'\in S, {\mathbf x}(\xi
),{\mathbf x}(\xi')\in U, q(\xi,\xi')>0}|\bar{\mathbf x}^i(\xi
)-\bar{\mathbf x}^i(\xi')|,\qquad i=1,\ldots,d,
\]
and use (\ref{CTB}) to see that $J_i\le\bar J\s_i$. Determine
$\theta_i\in(0,\infty)$ by $\theta_i e^{\theta_iJ_i}={\delta }/(\rho
t_0\s_i)$; then $\theta_i\le{\delta}/(\rho t_0\s_i)$, so
$\theta_iJ_i\le2{\delta }\bar J/(3 \bar A t_0)\le1/4$, since
we\vspace*{1pt} assumed that $\delta\bar J \le\bar A t_0/4$. Since
$e^{1/4}\le4/3$, we have $\rho e^{\theta _iJ_i}\le2\bar A$. We now
apply the exponential martingale inequality, Proposition~\ref{EMI},
substituting $\pm\bar{\mathbf x}^i$ for $\phi$ for $i=1,\ldots,d$ and
substituting ${\delta}\s_i$ for ${\delta}$ and $\rho t_0\s_i^2$ for
$\ve$. We thus obtain
\[
\PP\bigl(\O(M)^c\cap\O(a)\bigr)\le2de^{-{\delta}^2/(4\bar At_0)}.
\]
If we take $\bar A=A$, then, using (\ref{ALA}) and (\ref{LBTT}), we
have $\O(a)=\O$, so the proof of (\ref{FLEA}) is now complete.

Set $\eta=16\L\t^2A^2$. We shall complete the proof of (\ref
{FLEAB}) by showing that
\[
\PP\bigl(\O(a)^c\cap\O(\g,g)\bigr)\le2de^{-\tilde{\delta}^2/(4\eta t_0)}.
\]
Take $\tilde M$ as in (\ref{KEYT}), with $a$ as in (\ref
{ALA}). Then,
for $t\le T$,
%
%
\begin{eqnarray}\label{KEYT-1}
\int_0^ta(Y_s)\,ds&=&\tilde\chi(\bX_t,Y_t)-\tilde\chi(\bX
_0,Y_0)-\tilde M_t\nonumber\\[-8pt]\\[-8pt]
&&{}+\int_0^t\tilde\D(X_s)\,ds
+\int_0^t\bar a(\bX_s)\,ds,\nonumber
\end{eqnarray}
where $\tilde\D=G\tilde\chi({\mathbf x},y)-Q(\tilde\chi({\mathbf x},y))$.
Consider the event
\[
\O(\tilde M)=\Bigl\{\sup_{t\le T_0\wedge t_0}|\tilde M_t^i|\le
\tilde{\delta}\s_i^2 \mbox{ for } i=1,\ldots,d\Bigr\}.
\]
Then, on $\O(\g,g)\cap\O(\tilde M)$, we can estimate the terms in
(\ref{KEYT-1}) to obtain
\[
\int_0^{T_0\wedge t_0}a^i(Y_s)\,ds\le(5\tilde{\delta}+\bar At_0)\s
_i^2\le\rho t_0\s_i^2.
\]
Hence, it will suffice to show that
\[
\PP(\O(\tilde M)^c)\le2de^{-\tilde{\delta}^2/(4\eta t_0)}.
\]
For this, we again use the exponential martingale inequality.
Take $\phi(\xi)=\pm\tilde\chi^i({\mathbf x}(\xi),y(\xi))$ in
Proposition~\ref{EMI} and note that $\a^\phi(\xi)\le16\L\t
^2A^2\s_i^4$, so
\[
\int_0^{T_0\wedge t_0}\a^\phi(X_s)\,ds\le16\L\t^2A^2\s_i^4t_0=\eta
t_0\s_i^4.
\]
Set
\[
\tilde J_i=J(\phi)=\sup_{\xi,\xi'\in S, {\mathbf x}(\xi),{\mathbf
x}(\xi')\in U, q(\xi,\xi')>0}|\phi(\xi)-\phi(\xi')|,\qquad
i=1,\ldots,d,
\]
then $\tilde J_i\le4\t A\s_i^2$.
Determine $\tilde\theta_i\in(0,\infty)$ by $\tilde\theta_i e^{\tilde
\theta_i \tilde J_i}=\tilde{\delta}/(\eta t_0\s_i^2)$.
Then $\tilde\theta_i\le\tilde{\delta}/(\eta t_0\s_i^2)$ so $\tilde
\theta_i \tilde J_i\le\bar A/(40\L\t A)\le1/2$ and so
$e^{\tilde\theta_i \tilde J_i}\le2$.
Hence,
\[
\PP(\O(\tilde M)^c)\le2d\exp\{-\tilde{\delta}^2/(2\eta
t_0e^{\tilde\theta_i \tilde J_i})\}\le2de^{-\tilde{\delta
}^2/(4\eta t_0)}
\]
as required.
\end{pf}

\section{The supermarket model with memory}\label{SMM}

The supermarket model with memory is a variant,
introduced in~\cite{MPS}, of the ``join the shortest queue'' model,
which has been widely studied~\cite{VDK,G,GA,LM,LMA,MR2152252}.
We shall rigorously verify the asymptotic picture for large numbers of
queues derived in~\cite{MPS}.
This will serve as an example to illustrate the general theory
of the preceding sections. The explicit form of the error probabilities
in Theorem~\ref{FLET}
is used to advantage in dealing with the infinite-dimensional character
of the limit model.

Fix $\lambda\in(0,1)$ and an integer $n\ge1$. We shall consider the limiting
behavior as $N\to\infty$ of the following queueing system.
Customers arrive as a Poisson process of rate $N \lambda$ at a system
of $N$ single-server queues. At any given time, the length of one of
the queues is kept
under observation. This queue is called the \textit{memory queue}.
On each arrival, an independent random sample of size $n$ is chosen
from the set of all
$N$ queues. For simplicity, we sample with replacement, allowing
repeats and allowing the choice of the memory queue.
The customer joins whichever of the memory queue or the sampled queues
is shortest, choosing randomly in the event of a tie.
Immediately after the customer has joined a queue, we switch the memory queue,
if necessary, so that it is the currently shortest queue among the
queues just sampled and the previous memory queue.
The service requirements of all customers are assumed independent and
exponentially distributed of mean $1$.

Write $Z_t^k=Z_t^{N,k}$ for the proportion of
queues having at least $k$ customers at time $t$, and write $Y_t$ for length
of the memory queue at time $t$. Set $Z_t=(Z^k_t\dvtx k \in\N)$ and
$X_t=(Z_t,Y_t)$.
Then $X=(X_t)_{t \ge0}$ is a Markov chain, taking values in
$S=S_0\times\Z^+$, where $S_0$ is the set of nonincreasing sequences
in $N^{-1}\{0,1,\ldots,N\}$ with finitely many nonzero terms.
We shall treat $Y$ as a fast variable and prove a fluid limit for $Z$
as $N\to\infty$.

\subsection{Statement of results}\label{SR}
Let $D$ be the set of nonincreasing sequences\footnote{To lighten the
notation, we shall sometimes move the coordinate index from a
superscript to a subscript, allowing the $n$th power of the $k$th
coordinate to be written $z_k^n$. We shall also write the time
variable sometimes as a subscript, sometimes as an argument.}
$z=(z_k\dvtx k\in\N)$ in the interval $[0,1]$ such that
\[
m(z):=\sum_kz_k<\infty.
\]
Define for $z\in D$ and $k\in\N$
%
%
\begin{equation}\label{MZK}
\mu(z,k)=\prod_{j=1}^k\frac{z_j^n}{1-p_{j-1}(z)},
\end{equation}
where
\[
p_{k-1}(z)=n(z_{k-1}-z_k)z_k^{n-1}
\]
and where we take $z_0=1$.
Set $\mu(z,0)=1$ for all $z$.
An elementary calculation (maximizing over $z_k$ while keeping
$z_{k-1}$ fixed) shows that in the case $n\ge2$,
%
%
\begin{equation}\label{PKZ}
p_{k-1}(z)\le z_{k-1}^n (1-1/n)^{n-1}\le(1-1/n)^{n/2} \le e^{-1/2}<1.
\end{equation}
In the case $n=1$ we have $p_{k-1}(z)=z_{k-1}-z_k\le1$ and it is
possible that $0/0$ appears in the product (\ref{MZK}).
For definiteness we agree to set $0/0=1$ in this case.
Note that $\mu(z,k)\ge\mu(z,k+1)$ for all $k\ge0$.
Define for $z\in D$
\[
v_k(z)={\lambda}z_{k-1}^n\mu(z,k-1)-{\lambda}z_k^n\mu(z,k)-(z_k-z_{k+1})
\]
and consider the differential equation
%
%
\begin{equation}\label{ODE}
\dot z(t)=v(z(t)),\qquad t\ge0.
\end{equation}
By a \textit{solution} to (\ref{ODE}) in $D$ we mean a family of
differentiable functions $z_k\dvtx[0,\infty)\to[0,1]$
such that for all $t\ge0$ and $k\in\N$ we have $(z_k(t)\dvtx k\in\N)\in
D$ and
\[
\dot z_k(t)=v_k(z(t)).
\]

%
\begin{theorem}\label{DDD}
For all $z(0)\in D$, the differential equation $\dot z(t)=v(z(t))$ has
a unique
solution in $D$ starting from $z(0)$. Moreover, if $(w(t)\dvtx t\in D)$ is
another solution
in $D$ with $z_k(0)\le w_k(0)$ for all $k$, then $z_k(t)\le w_k(t)$ for
all $k$ and all $t\ge0$.
\end{theorem}

There is a fixed point of these dynamics $a\in D$ given by setting
$a_0=1$ and defining
%
%
\begin{equation}\label{DEFA}
a_{k+1}={\lambda}a_k^n\mu(a,k),\qquad k\ge0.
\end{equation}
The components of $a$ decay super-geometrically.
Set
%
%
\begin{equation}\label{ALP}
\a=n+\tfrac12+\sqrt{n^2+\tfrac14}.
\end{equation}
Then $\a\in(2n,2n+1)$.
%
%
\begin{theorem}\label{SGD}
We have
\[
\lim_{k\to\infty}\frac1k\log\log\biggl(\frac1{a_k}\biggr)=\a.
\]
\end{theorem}

Assume for simplicity that we start the queueing system from the state
where all queues, except the
memory queue, are empty and where the memory queue has exactly one
customer. Write
$(z(t)\dvtx t\ge0)$ for the solution to (\ref{ODE}) starting from~$0$.
Then $z_k(t)\le a_k$ for all
$k$ and $t$. Our main result shows that $(z(t)\dvtx t\ge0)$ is a good
approximation to the process of empirical
distributions of queue lengths $(Z^N(t)\dvtx t\ge0)$ for large $N$. The
sense of this approximation is reasonably sharp.
In particular, as a straightforward corollary, we obtain that, on a
given time interval $[0,t_0]$,
for any $r>\a^{-1}$,
with high probability as $N\to\infty$, no queue length exceeds $r\log
\log N$.
%
%
\begin{theorem}\label{MR}
Set $\k=(2\a)^{-1}$ and define
\[
d=d(N)=\sup\{k\in\N\dvtx Na_k>N^{\kappa}\}.
\]
Fix a function $\phi$ on $\N$ such that $\phi(N)/N^\k\to0$ and
$\log\phi(N)/\log\log N\to\infty$ as
$N\to\infty$.
Set $\rho=4/(1-{\lambda})$ when $n=1$ and set $\rho
=2^n/(1-e^{-1/2})$ when $n\ge2$.
Set $\tilde a_{d+1}=N^{-1}a_d^n+\rho^da_{d+1}$.
Then
%
%
\begin{equation}\label{DLLN}
\lim_{N\to\infty}d(N)/\log\log N=1/\a.
\end{equation}
Moreover, for all $t_0\ge0$, we have
%
%
\begin{equation}\label{MRA}
\lim_{N\to\infty}\PP\Biggl(\sup_{t\le t_0}\sup_{k\le d}\frac
{|Z_t^{N,k}-z_t^k|}{\sqrt{a_k}}\ge\sqrt{\frac{\phi(N)}{N}}\Biggr)=0
\end{equation}
and
%
%
\begin{equation}\label{MRB}
\lim_{R\to\infty}\limsup_{N\to\infty}\PP(Z_t^{N,d+1}\ge R\tilde
a_{d+1}\mbox{ for some }t\le t_0)=0
\end{equation}
and
%
%
\begin{equation}\label{MRC}
\lim_{N\to\infty}\PP(Z^{N,d+2}_t=0\mbox{ for all }t\le t_0)=1.
\end{equation}
\end{theorem}

The argument used to prove this result would apply without modification
starting from any initial condition $z(0)$
for the limit dynamics (\ref{ODE}) such that $z_k(0)\le a_k$ for all
$k$, with suitable conditions
on the convergence of $Z^N(0)$ to $z(0)$.
It may be harder to move beyond initial conditions which do not lie
below the fixed point. We do note here, however,
a family of long-time upper bounds for the limit dynamics which might
prove useful for such an extension.
Fix $j\in\N$ and define
$a^{(j)}_k=a_{(k-j)^+}$ for each $k \in\Z^+$; then $a^{(j)}$ is a
fixed point of the modified equation
\[
\dot w_k(t)=v_k(w(t))+\bigl(w_j(t)-w_{j+1}(t)\bigr)1_{\{k=j\}}.
\]
Since the added term is always nonnegative, a similar argument to that
used to prove Theorem~\ref{DDD} in the next subsection
also shows that, if $z(0)\le a^{(j)}$ and $(z(t)\dvtx t\ge0)$ is a solution
of the original equation, then $z(t)\le a^{(j)}$ for all~$t$.

\subsection{Existence and monotonicity of the limit dynamics}\label{EMLD}
The differential equation (\ref{ODE}) characterizes the limit dynamics
for the fluid variables in our queueing model.
Our analysis of its space of solutions will rest on the exploitation
of certain nonnegativity properties which have a natural probabilistic
interpretation.
We shall use the following standard property of differential equations:
if $b=(b^1,\ldots,b^d)$ is a Lipschitz vector field
on $\R^d$ such that $b^1(x)\ge0$ whenever $x=(x^1,\ldots,x^d)$ with
$x^1=0$, and if $(x(t)\dvtx t\le t_0)$ is a solution
to $\dot x(t)=b(x(t))$ with $x^1(0)\ge0$, then $x^1(t)\ge0$ for all
$t\le t_0$.\vadjust{\goodbreak}

We consider first a truncated, finite-dimensional system.
Fix $d\in\N$ and define a vector field $u=u^{(d)}$ on $D$ by setting
$u_k(z)=v_k(z)$ for $k\le d-1$ and
%
%
\begin{equation}\label{UDZ}
u_d(z)={\lambda}z_{d-1}^n\mu(z,d-1)-{\lambda}z_d^n\mu(z,d)-z_d
\end{equation}
and $u_k(z)=0$ for $k\ge d+1$.
Set $D(d)=\{(x_1,\ldots,x_d,0,0,\ldots)\dvtx0\le x_d\le\cdots\le x_1\le
1\}$.
%
%
\begin{proposition}\label{ADD}
For all $x(0)\in D(d)$, the differential equation $\dot x(t)=u(x(t))$
has a unique
solution $(x(t)\dvtx t\ge0)$ in $D(d)$ starting from $x(0)$.
\end{proposition}
\begin{pf}
In the proof, we consider $D(d)$ as a subset of $\R^d$. The function
$u$ is continuous on $D(d)$ and is
differentiable in the interior of $D(d)$ with bounded partial derivatives.
[In the case $n=1$, the singularity in $(\pd/\pd x_j)\mu(x,k)$ for
$j\le k$ as $x_{j-1}-x_j\to1$ is
canceled by the factor $x_k$ by which it is multiplied, since $x_k\le
x_j$ on $D(d)$.]
For $x\in D(d)$ we have $u_1(x)\le0$ when $x_1=1$, and $u_d(x)\ge0$
when $x_d=0$.
Moreover, for $k=1,\ldots,d-1$, if $x_k=x_{k+1}$ then $p_k(x)=0$ so
\begin{eqnarray*}
u_{k+1}(x)&=&x_k^n\bigl(\mu(x,k)-\mu(x,k+1)\bigr)
\le \mu(x,k)-\mu(x,k+1)\\
&\le& x_{k-1}^n\mu(x,k-1)-x_k^n\mu(x,k)\le u_k(x).
\end{eqnarray*}
The conclusion now follows by standard arguments.
\end{pf}
\begin{pf*}{Proof of Theorem~\ref{DDD}}
Suppose that $(w(t)\dvtx t\ge0)$ is a solution to $\dot w(t)=v(w(t))$
in $D$ starting\vspace*{1pt} from $w(0)$, with $z(0)\le w(0)$, that is
to say $z_k(0)\le w_k(0)$ for all $k$. Fix $d$ and write
$x(t)=z^{(d)}(t)$ for the solution to $\dot x(t)=u^{(d)}(x(t))$ in
$D(d)$ starting from $(z_1(0),\ldots,z_d(0),0,0,\ldots)$. Set
$y(t)=(w_1(t),\ldots,w_d(t),0,0,\ldots)$ and note that $x(0)\le y(0)$
and $y(t)\in D(d)$ for all $t$. We shall show that $x(t)\le y(t)$ for
all $t$. Now consider $D(d)$ as a subset of $\R^d$. We have
\[
\dot y(t)=u(y(t))+w^{d+1}(t)e_d,
\]
where $e_d=(0,\ldots,0,1)$. Note that $w^{d+1}(t)\ge0$ for all $t$.
Now $u$ is Lipschitz on $D(d)$ and for $k=1,\ldots,d$ we can show
that\footnote{An elementary calculation shows that
$\mu(x,j)\le\mu(y,j)$ for all $j$ whenever $x\le y$.
This will also be shown by a soft probabilistic argument in Section
\ref{LED}.
The further condition $x_k=y_k$ gives the inequality
\[
y_{k-1}^n-x_{k-1}^n\ge n(y_{k-1}-x_{k-1})x_k^{n-1}=p_{k-1}(y)-p_{k-1}(x)
\ge\frac{y_k^{2n}}{1-p_{k-1}(y)}-\frac{x_k^{2n}}{1-p_{k-1}(x)}.
\]
}
\[
x,y\in D(d),\qquad x\le y,\qquad
x_k=y_k \quad\implies\quad u_k(x)\le u_k(y).
\]
Hence, by a standard argument $z^{(d)}(t)=x(t)\le y(t)\le w(t)$ for all $t$.
The same argument shows that $z^{(d)}(t)\le z^{(d+1)}(t)$ for all $t$,
so the limit
$z_k(t)=\lim_{d\to\infty}z^{(d)}_k(t)$ exists for all $k$ and $t$,
and $z(t)\le w(t)$ for all $t$.\vadjust{\goodbreak}

Fix $k$ and take $d\ge k+1$. Then the following equation holds for all $t$:
\begin{eqnarray*}
&&z^{(d)}_k(t)+\int_0^t{\lambda}z_k^{(d)}(s)^n\mu
\bigl(z^{(d)}(s),k\bigr)\,ds+\int_0^tz_k^{(d)}(s)\,ds\\
&&\qquad =z_k(0)+\int_0^t{\lambda}z_{k-1}^{(d)}(s)^n\mu
\bigl(z^{(d)}(s),k-1\bigr)\,ds+\int_0^tz_{k+1}^{(d)}(s)\,ds.
\end{eqnarray*}
On letting $d\to\infty$, we see by monotone convergence that
\begin{eqnarray*}
&&z_k(t)+\int_0^t{\lambda}z_k(s)^n\mu(z(s),k)\,ds+\int_0^tz_k(s)\,ds\\
&&\qquad =z_k(0)+\int_0^t{\lambda}z_{k-1}(s)^n\mu\bigl(z(s),k-1\bigr)\,ds+\int
_0^tz_{k+1}(s)\,ds.
\end{eqnarray*}
Since $z(t)\in D$ for all $t$, all integrands in this equation are
bounded by $1$. It is now
straightforward to see that $(z(t)\dvtx t\ge0)$ is a solution.

Now
\begin{eqnarray*}
&&w_k(t)+\int_0^t{\lambda}w_k(s)^n\mu(w(s),k)\,ds+\int_0^tw_k(s)\,ds\\
&&\qquad =w_k(0)+\int_0^t{\lambda}w_{k-1}(s)^n\mu\bigl(w(s),k-1\bigr)\,ds+\int
_0^tw_{k+1}(s)\,ds.
\end{eqnarray*}
By summing these equations over $k\in\{1,\ldots,d\}$ we see that the
map $t\mapsto\sum_{k=1}^d w_k(t)-{\lambda}t$
is nonincreasing for all $d$. Hence, $m(w(t))\le m(w(0))+{\lambda
}t<\infty$. The equations can then be summed
over all $k$ and rearranged to obtain
\[
m(w(t))=m(w(0))+{\lambda}t-\int_0^tw_1(s)\,ds.
\]
On the other hand,
\[
m\bigl(z^{(d)}(t)\bigr)=m\bigl(z^{(d)}(0)\bigr)+{\lambda}t-\int
_0^tz_1^{(d)}(s)\,ds-{\lambda}\int_0^tz^{(d)}_d(s)^n\mu\bigl(z^{(d)}(s),d\bigr)\,ds
\]
so
%
%
\begin{eqnarray}\label{MDM}
m(w(t))-m\bigl(z^{(d)}(t)\bigr)
&\le& m(w(0))-m\bigl(z^{(d)}(0)\bigr)\nonumber\\[-8pt]\\[-8pt]
&&{}+{\lambda}\int_0^t
z_d(s)^n\mu(z(s),d)\,ds.\nonumber
\end{eqnarray}
If $w(0)=z(0)$ then the right-hand side tends to $0$ as $d\to\infty$
so we must have $z(t)=w(t)$ for all $t$.
\end{pf*}

\subsection{Properties of the fixed point}
Recall the definition (\ref{DEFA}) of the fixed point~$a$. Since $\mu
(a,k)\le1$ for all $k$, we have
\[
a_k\le{\lambda}^{1+n+\cdots+n^{k-1}}\vadjust{\goodbreak}
\]
so $a_k\to0$ as $k\to\infty$.
Theorem~\ref{SGD} is then a straightforward corollary of the following
estimate.
%
%
\begin{proposition}\label{AK}
There is a constant $C({\lambda},n)<\infty$ such that, for all $k\ge0$,
%
%
\begin{equation}
C^{-1}a_k^\a\le a_{k+1}\le Ca_k^\a,
\end{equation}
where $\a$ is given by (\ref{ALP}).
\end{proposition}
\begin{pf}
Note that, since $a_1={\lambda}$, we have $p_{k-1}(a)=a_{k-1}-a_k\le
{\lambda}\vee(1-{\lambda})<1$ for all $k$ when $n=1$. On the other
hand, equation (\ref{PKZ}) gives $p_{k-1}(z)\le e^{-1/2}$ for all $k$
when $n\ge2$. Then from $\sum_{k=1}^\infty p_{k-1}(a)\le1$,
we obtain a constant $c<\infty$, which may depend on ${\lambda}$ when
$n=1$, such that
\[
\prod_{k=1}^\infty\frac1{1-p_{k-1}(a)}\le c.
\]
Then for $k\ge0$
\[
{\lambda}a_k^{2n}\prod_{j=1}^{k-1}a_j^n\le a_{k+1}\le c{\lambda
}a_k^{2n}\prod_{j=1}^{k-1}a_j^n,
\]
so for $k\ge1$,
%
%
\begin{equation}\label{CAK}
c^{-1}a_k^{2n+1}a_{k-1}^{-n}\le a_{k+1}\le ca_k^{2n+1}a_{k-1}^{-n}.
\end{equation}
Note that ${\lambda}a_0^\a={\lambda}=a_1\le{\lambda}^{-1}a_0^\a$.
Fix $A\ge1/{\lambda}$ and suppose inductively that
\[
A^{-1}a_{k-1}^\a\le a_k\le Aa_{k-1}^\a.
\]
On using these inequalities to estimate $a_{k-1}$ in (\ref{CAK}), we obtain
\[
(cA^{n/\a})^{-1}a_k^\a\le a_{k+1}\le cA^{n/\a}a_k^\a,
\]
where we have used the fact that $1-n/\a=\a-2n$. Hence, the induction
proceeds provided
we take $A\ge c^{\a/(\a-n)}$.
\end{pf}

\subsection{Choice of fluid coordinates and fast variable}
In the remaining subsections we apply Theorem~\ref{FLET} to deduce
Theorem~\ref{MR}.
Define $d$ as in Theorem~\ref{MR} and take as auxiliary space $I=\N$
when $n=1$ and $I=\Z^+$ when $n\ge2$.
Make the following choice of fluid and auxiliary coordinates:
for $\xi=(z,y)\in S$ with $z=(z_k\dvtx k\in\N)$, set
\[
x^k(\xi)=z_k,\qquad k=1,\ldots,d,\qquad y(\xi)=y.
\]
Thus our fluid variable is $\bX_t={\mathbf x}(X_t)=(Z_t^1,\ldots
,Z_t^d)$ and our fast variable is $Y_t=y(X_t)$.
Note that when $n=1$, if $Y_0\ge1$ then $Y_t\ge1$ for all $t$, so $Y$
takes values in $I=\N$.

Let us compute the drift vector $\b(\xi)$ for $\bX$ when $X$ is in
state $\xi=(z,y)\in S$.
Note that $X^k$ makes a jump of size $1/N$ when a customer arrives at a
queue of length $k-1$,
and makes a jump of size $-1/N$ when a customer departs from a queue of
length $k$,
otherwise $X^k$ is constant.
The length of the queue which an arriving customer joins depends on the
length of the memory queue $y$ and on the
lengths of the sampled queues. Denote the vector of sampled queue
lengths by $V=V(z)=(V_1,\ldots,V_n)$
and write $V^{(1)}\le V^{(2)}\le\cdots\le V^{(n)}$ for the ordered
queue lengths.
Define $\min(v)=v_1\wedge\cdots\wedge v_n$ and set $M=\min(V)$. Then
$M=V^{(1)}$ and
\[
\PP(M\ge k)=z_k^n.
\]
A new customer will go to a queue of length at least $k$ if and only if
$M\ge k$ and $y\ge k$.
So the rate for an arrival to a queue of length exactly $k-1$ is
\[
N{\lambda}\PP(M\ge k-1)1_{\{y\ge k-1\}}-N{\lambda}\PP(M\ge k)1_{\{
y\ge k\}}.
\]
The rate for a departure from a queue of length $k$ is $N(z_k-z_{k+1})$.
Hence, setting $z_0=1$, we have
\[
\b_k(\xi)={\lambda}z_{k-1}^n1_{\{y\ge k-1\}}-{\lambda}z_k^n1_{\{
y\ge k\}}-(z_k-z_{k+1}).
\]

We now compute (an approximation to) the jump rates $\g(\xi,y')$ for
$Y$ when $X$ is in state $\xi=(z,y)\in S$.
The rate of departures from the memory queue is at most~$1$.
Arrivals to the system occur at rate $N{\lambda}$. Occasionally, the
memory queue falls in the sample,
an event of probability no greater than $n/N$ and hence, of rate no
greater than ${\lambda}n$.
Assuming that the the memory queue does not fall in the sample,
the length of the memory queue after an arrival is given by
%
%
\begin{equation}\label{YP}
F(y,V)=(y+1)1_{\{y\le M-1\}}+y1_{\{M\le y\le P\}}+P1_{\{y\ge P+1\}},
\end{equation}
where $P=p(V)$ is given by $P=M+1$ when $n=1$ and $P=(M+1)\wedge
V^{(2)}$ otherwise.
Hence, we have
%
%
\begin{equation}\label{GAM}
\sum_{y'\not=y(\xi)}\bigl|\g(\xi,y')-N{\lambda}\PP
\bigl(F(y,V(z))=y'\bigr)\bigr|\le1+{\lambda}n.
\end{equation}

\subsection{Choice of limit characteristics and coupling
mechanism}\label{CLC}
Define
\[
U=\{x\in\R^d\dvtx0\le x_d\le\cdots\le x_1\le1\mbox{ and }x_1\le
({\lambda}+1)/2\mbox{ and }x_k\le2a_k\mbox{ for all }k\}.
\]
The condition $x_1\le({\lambda}+1)/2$ ensures that $1-x_1$ is
uniformly positive on $U$.
Define $b\dvtx U\times\Z^+\to\R^d$ by
%
%
\begin{equation}\label{BKX}
b_k(x,y)={\lambda}x_{k-1}^n1_{\{y\ge k-1\}}-{\lambda}x_k^n1_{\{y\ge
k\}}-(x_k-x_{k+1}),
\end{equation}
where we set $x_0=1$ and $x_{d+1}=0$.
Then, for $\xi\in S$ with ${\mathbf x}(\xi)\in U$, we have
%
%
\begin{equation}\label{BETAB}
\b(\xi)=b({\mathbf x}(\xi),y(\xi))+(0,\ldots,0,z_{d+1}).
\end{equation}
It is convenient to specify our choice of the generator matrices
$(G_x\dvtx x\in U)$ and our choice of
coupling mechanism at the same time.
Set $\nu=N{\lambda}$ and take as auxiliary space $E=(\Z^+)^n$.
Define a family of probability
distributions $\mu=(\mu_x\dvtx x\in U)$ on~$E$, taking $\mu_x$ to be the
law of a random sample
$V=V(x)=(V_1,\ldots,V_n)$
with
\[
\PP(V_1\ge k)=\cdots=\PP(V_n\ge k)=x_k
\]
for $k=0,1,\ldots,d+1$. Note that
%
%
\begin{equation}\label{MTV}
\|\mu_x-\mu_{x'}\|\le2n\sum_{k=1}^d|x_k-x'_k|.
\end{equation}
Then define for distinct $y,y'\in\Z^+$,
\[
g(x,y,y')=N{\lambda}\PP\bigl(F(y,V(x))=y'\bigr),
\]
where $F$ is given by (\ref{YP}). We take as coupling mechanism the
triple $(\nu,\mu,F)$.

Note that $F(y,v)=F(\bar y,v)$ for all $y,\bar y\in I$ whenever
$p(v)=\min I$. For $x\in U$ we have
\[
\PP\bigl(p(V(x))=1\bigr)\ge1-x_1>\frac{1-{\lambda}}2,
\]
when $n=1$, whereas for $n\ge2$ we have
\[
\PP\bigl(p(V(x))=0\bigr)\ge(1-x_1)^2>\biggl(\frac{1-{\lambda
}}2\biggr)^2.
\]
Hence, we obtain, in all cases, $m(x,y,\bar y)\le\t$, where we set
\[
\t=\frac4{N{\lambda}(1-{\lambda})^2}.
\]

For $\xi\in S$ with $x={\mathbf x}(\xi)\in U$ we can realize a sample
$V(z)$ (from the distribution of queue lengths)
and the sample $V(x)$ on the same probability space by setting
$V_i(x)=V_i(z)\wedge d$.
Write $M(x)=\min(V(x))$ and $P(x)=p(V(x))$.
Then $M(x)=M(z)\wedge d$ and
$P(x)=P(z)\wedge(d+1)$ when $n=1$ and $P(x)=P(z)\wedge d$ when $n\ge2$.
The difference between the two cases is that there is no second
shortest queue in the sample when $n=1$.
We have, for $n=1$,
\[
\PP\bigl(P(z)\not= P(x)\bigr)\le\PP\bigl(M(z)\ge d+1\bigr)=z_{d+1}
\]
and, for $n\ge2$,
\[
\PP\bigl(P(z)\not= P(x)\bigr)\le\PP\bigl(P(z)\ge d+1\bigr)\le nz_dz_{d+1}^{n-1}.
\]
Now $P(x)=P(z)$ implies $M(x)=M(z)$ and hence, $F(y,V(x))=F(y,V(z))$
for all $y$. Hence,
\[
\PP\bigl(F(y,V(z))\not=F(y,V(x))\bigr)\le\PP\bigl(P(x)\not=P(z)\bigr).
\]
On combining this with (\ref{GAM}) we obtain
%
%
\begin{equation}\label{GAMMAG}
\sum_{y'\not=y(\xi)}|\g(\xi,y')-g({\mathbf x}(\xi),y(\xi
),y')|\le1+{\lambda}n+N{\lambda}nz_{d+1}.
\end{equation}

\subsection{Local equilibrium distribution}\label{LED}
The Markov chain determined by the generator $G_x$ has a unique closed
communicating class, which is contained in $\{0,1,\ldots,d\}$.
Hence, $G_x$ has a unique equilibrium
distribution $\pi_x$ which is supported on $\{0,1,\ldots,d\}$.
Consider a continuous-time Markov chain\footnote{See footnote \ref
{NFV}.} $Y=(Y_t)_{t\ge0}$ with generator $G_x$, and initial
distribution $\pi_x$. Set $\mu(x,k)=\PP(Y_0\ge k)$.
Then $Y$ jumps into $\{0,1,\ldots,k\}$ from $j$ at rate
$\alpha=N{\lambda}(1-x_{k+1}^n-p_k(x))$, for all $j\ge k+1$. On the
other hand, $Y$ jumps out of
$\{0,1,\ldots,k\}$ only from $k$, and that at rate $\beta=N{\lambda
}x_{k+1}^n$. Since the long run rates
of such jumps must agree, we deduce that $\alpha\mu(x,k+1)=\beta\pi_x(k)$.
Hence, we obtain
\[
\mu(x,k+1)\bigl(1-p_k(x)\bigr)=x_{k+1}^n\mu(x,k)
\]
and so
%
%
\begin{equation}\label{IMT}
\mu(x,k)=\prod_{j=1}^k\frac{x_j^n}{1-p_{j-1}(x)},\qquad k=1,\ldots,d.
\end{equation}
Hence, our present notation is consistent with the definition (\ref{MZK}).
Note also that, for $z\in D$, $\mu(z,k)$ depends only on $z_1,\ldots
,z_k$; in particular,
if $x=(z_1,\ldots,z_d)$ then $\mu(z,k)=\mu(x,k)$ for all $k\le d$.
Note that $\bar b$ is given by
%
%
\begin{equation}\label{BARBE}
\bar b_k(x)={\lambda}x_{k-1}^n\mu(x,k-1)-{\lambda}x_k^n\mu
(x,k)-(x_k-x_{k+1}),\qquad k=1,\ldots,d,\hspace*{-35pt}
\end{equation}
where $x_0=1$ and $x_{d+1}=0$. Hence, $\bar b=u^{(d)}$ as defined in
Section~\ref{EMLD}.

A comparison of (\ref{UDZ}) and (\ref{BKX}) now shows that $\bar b=u^{(d)}$.

Recall that $\rho=4/(1-{\lambda})$ when $n=1$ and that $\rho
=2^n/(1-e^{-1/2})$ when $n\ge2$.
Then for $x\in U$ and $k\ge1$ we have
\begin{eqnarray*}
\mu(x,k)&=&x_k^n\mu(x,k-1)/\bigl(1-p_{k-1}(x)\bigr)\le\rho a_k^n\mu(x,k-1)\\
&\le&
\rho a_k^n\mu(x,k-1)/\bigl(1-p_{k-1}(a)\bigr)
\end{eqnarray*}
so, for all $x\in U$ and inductively for all $k\ge1$, we obtain
%
%
\begin{equation}\label{MUXA}
\mu(x,k)\le\rho^k\mu(a,k).
\end{equation}

The following argument shows that $\mu(z,k)\le\mu(z',k)$ for all $k$
whenever $z\le z'$.
Fix $d\ge k$ and set $x'=(z_1',\ldots,z_d')$. Assume that $x\le x'$. By
a standard construction
we can realize samples $V=V(x)$ and $V'=V(x')$ on a common probability
space such that $V_i\le V'_i$
for all $i$. Then we can construct Markov chains $Y$ and $Y'$,
having generators $G_x$ and $G_{x'}$, respectively, on the canonical
space of a marked Poisson process
of rate $N{\lambda}$, where the marks are independent copies of
$(V,V')$, as follows.
Set $Y_0=Y'_0=1$ and define recursively at each jump time $T$ of the
Poisson process
$Y_T=F(Y_{T-},V_T)$ and $Y_T'=F(Y'_{T-},V_T')$, where $(V_T,V_T')$ is
the mark at time~$T$.
Then, since $F$ is nondecreasing in both arguments, we see by induction
that $Y_t\le Y'_t$ for all $t$.
Hence, by convergence to equilibrium,
\[
\mu(z,k)=\mu(x,k)=\lim_{t\to\infty}\PP(Y_t\ge k)\le\lim_{t\to
\infty}\PP(Y_t'\ge k)=\mu(x',k)=\mu(z',k).
\]

\subsection{Corrector upper bound}
We take as our reference state $\bar y=\min I$ and note that, under the
coupling mechanism, we have $\bar Y_t\le Y_t$ for all $t$.
Then the $k$th component of the corrector for $b$ is given by
\[
\chi_k(x,y)={\lambda}\E_y\int_0^{T_c}\bigl(x_{k-1}^n1_{\{\bar Y_s<k-1\le
Y_s\}}-x_k^n1_{\{\bar Y_s<k\le Y_s\}}\bigr)\,ds,
\]
so for $x\in U$ and all $y\in I$ we have
\[
|\chi_k(x,y)|\le\t x_{k-1}^n\le Ca_{k-1}^n/N.
\]
Now fix $y\le k-2$ and consider the stopping time $T=\inf\{t\ge
0\dvtx Y_t=k-1\}$. Note that $Y$ can enter state
$k-1$ only from $k-2$ and does so at rate $N{\lambda}x_{k-1}^n$,
whereas $T_c$ occurs in state $k-2$ at rate at least $N{\lambda
}(1-{\lambda})^2/4$.
Hence, $\PP_y(T\le T_c)\le Ca_{k-1}^n$ and so, by the Markov property,
\[
\E_y\int_0^{T_c}1_{\{\bar Y_s<k-1\le Y_s\}}\,ds\le\E_y\bigl(1_{\{T\le T_c\}
}m(x,k-1,\bar Y_T)\bigr)\le Ca_{k-1}^n\t\le C a_{k-1}^n/N.
\]
Hence, we obtain, for $x\in U$ and all $y\in I$,
%
%
\begin{equation}\label{CHID}
|\chi_k(x,y)|\le C\bigl(a_{k-1}^n1_{\{y\ge k-1\}}+a_{k-1}^{2n}\bigr)/N.
\end{equation}

\subsection{Quadratic variation upper bound}\label{QVB}
The growth rate at $\xi$ of the quadratic variation of the corrected
$k$th coordinate is given by
\[
\a_k(\xi)=\sum_{\xi'\not=\xi}\{\bar{\mathbf x}_k(\xi')-\bar
{\mathbf x}_k(\xi)\}^2q(\xi,\xi').
\]
Recall that $\bar x_k=x_k-\chi_k({\mathbf x},y)$. We estimate
separately, writing $x={\mathbf x}(\xi)$ and $y=y(\xi)$,
\[
\sum_{\xi'\not=\xi}\{{\mathbf x}_k(\xi')-x_k\}^2q(\xi,\xi')\le
N^{-2}\bigl(Nx_k+N{\lambda}x_{k-1}^n1_{\{y\ge k-1\}}\bigr)
\]
and
\[
\sum_{\xi'\not=\xi}\{\chi_k({\mathbf x}(\xi'),y(\xi'))-\chi
_k(x,y)\}^2q(\xi,\xi')\le CN^{-1}x_{k-1}^{2n}.
\]
When $y\le k-2$, we can improve the last estimate by splitting the sum
in two and using
\[
\sum_{\xi'\not=\xi, y(\xi')\le k-2}\{\chi_k({\mathbf x}(\xi
'),y(\xi'))-\chi_k(x,y)\}^2q(\xi,\xi')\le CN^{-1}x_{k-1}^{4n}
\]
and
\[
\sum_{\xi'\not=\xi, y(\xi')\ge k-1}\{\chi_k({\mathbf x}(\xi
'),y(\xi'))-\chi_k(x,y)\}^2q(\xi,\xi')\le CN^{-1}x_{k-1}^{3n}.
\]
We used (\ref{CHID}) for the first inequality and for the second used
\[
\sum_{\xi'\not=\xi, y(\xi')\ge k-1}q(\xi,\xi')\le CNx_{k-1}^n.
\]
On combining these estimates we obtain
\[
\a_k(\xi)\le C\bigl(x_k+x_{k-1}^n1_{\{y\ge k-1\}}+x_{k-1}^{3n}\bigr)\le
a_k(y(\xi))/N,
\]
where
\[
a_k(y(\xi))=C\bigl(a_k+a_{k-1}^n1_{\{y\ge k-1\}}\bigr)/N.
\]
Then, using the estimate (\ref{MUXA}) and the limit (\ref{DLLN}), we have
\[
\bar a_k(x)=C\bigl(a_k+a_{k-1}^n\mu(x,k-1)\bigr)\le C\rho^{d-1}a_k/N\le C(\log
N)^Ca_k/N,
\]
so we have for all $x\in U$ and $y\in I$
\[
a_k(y)\le Aa_k, \bar a_k(x)\le\bar Aa_k
\]
with $A=C/(Na_d^{1-n/\a})$ and $\bar A=C(\log N)^C/N$.
It is straightforward to check that, for $N$ sufficiently large, we
have $\bar A\le A\le\L\bar J^2$.

\subsection{Truncation estimates}
A specific feature of the problem we consider is that the limit
dynamics is infinite dimensional,
while the general fluid limit estimate applies in a finite-dimensional
context. In this subsection
we establish some truncation estimates which will allow us to reduce to
finitely many dimensions.

Let $(z(t)\dvtx t\ge0)$ be the solution in $D$ to $\dot z(t)=v(z(t))$
starting from $0$, as in
Theorem~\ref{MR}. Let $(x(t)\dvtx t\ge0)$ be the solution to $\dot
x(t)=\bar b(x(t))$ starting from $0$.
%
%
\begin{lemma}\label{TLD}
We have
\[
\sum_{k=1}^d|z_k(t)-x_k(t)|\le ta_{d+1}.
\]
\end{lemma}
\begin{pf}
Since $\bar b=u^{(d)}$, we have $x(t)=z^{(d)}(t)$ for all $t$, and so,
from (\ref{MDM}), we obtain
\begin{eqnarray*}
\sum_{k=1}^d|z_k(t)-x_k(t)|&\le&\sum_{k=1}^d\bigl(z_k(t)-z^{(d)}_k(t)\bigr)\le
{\lambda}\int_0^tz_d(s)^n\mu(z(t),d)\,ds\\
&\le& t{\lambda}a_d^n\mu(a,d)=ta_{d+1}.
\end{eqnarray*}
\upqed\end{pf}

Denote by $A_k(t)$ the number of arrivals to queues of length at least
$k$ by time $t$.
Note that $NZ^{k+1}_t\le A_k(t)$ for all $k\ge1$ and all $t$. Recall that
\[
\tilde a_{d+1}=N^{-1}a_d^n+\rho^da_{d+1}.
\]

%
\begin{lemma}\label{ADT}
There is a constant $C({\lambda},n)<\infty$ such that, for all $t\ge
0$ and all~$N$, we have
\[
\E\bigl(A_d(T\wedge t)\bigr)\le Ce^{Ct}N\tilde a_{d+1}.
\]
\end{lemma}
\begin{pf}
Consider the function $f$ on $U\times I$ given by $f(x,y)=1_{\{y\ge d\}
}$ and note that $\bar f(x)=\mu(x,d)$.
Let $\chi$ be the corrector for $f$ given by (\ref{CFF}).
Then, for all $x\in U$ and all $y\in I$,
\[
|\chi(x,y)|\le2\t\|f\|_\infty=2\t=CN^{-1}
\]
and, whenever $x={\mathbf x}(\xi)$ and $x'={\mathbf x}(\xi')$ with
$q(\xi,\xi')>0$, by the estimates (\ref{CEE}) and (\ref{MTV}),
%
%
\begin{equation}\label{DRT}
|\chi(x,y)-\chi(x',y)|\le2\nu\t^2\|f\|_\infty\|\mu_x-\mu_{x'}\|
\le CN^{-2}.
\end{equation}
By optional stopping,
\[
\biggl|\int_0^{T\wedge t}Q(\chi({\mathbf x},y))(X_s)\,ds\biggr|=\bigl|\E
\bigl(\chi(\bX_{T\wedge t},Y_{T\wedge t})-\chi(\bX_0,Y_0)\bigr)\bigr|
\le CN^{-1}.
\]
Now
\[
Q(\chi({\mathbf x},y))(\xi)=1_{\{y\ge d\}}-\mu({\mathbf x}(\xi
),d)-\D_1(\xi)-\D_2(\xi),
\]
where $\D_1,\D_2$ are given by (\ref{DGH}), (\ref{DGH-1}). We use
(\ref{GAMMAG}) to obtain the estimate
\[
|\D_1(\xi)|\le2\t(1+{\lambda}n+N{\lambda}nz_{d+1})=C(z_{d+1}+N^{-1})
\]
and from (\ref{DRT}) deduce that
\[
|\D_2(\xi)|\le N(1+{\lambda})CN^{-2}=CN^{-1}.
\]
So
\[
\E\int_0^{T\wedge t}1_{\{Y_s\ge d\}}\,ds
\le CN^{-1}+\E\int_0^{T\wedge t}\bigl(\mu(\bX
_s,d)+CZ_s^{d+1}+CN^{-1}\bigr)\,ds.
\]
Set $g(t)=\E(A_d(T\wedge t))$, then
\begin{eqnarray*}
g(t)&=&N{\lambda}\E\int_0^{T\wedge t}(X^d_s)^n1_{\{Y_s\ge d\}}\,ds\le
N{\lambda}2^na_d^n\E\int_0^{T\wedge t}1_{\{Y_s\ge d\}}\,ds\\
&\le& Ca_d^n\biggl(1+\int_0^t\bigl(N\rho^d\mu(a,d)+1+g(s)\bigr)\,ds\biggr)\\
&\le&
CN\tilde a_{d+1}(1+t)+C\int_0^tg(s)\,ds.
\end{eqnarray*}
Here we have used the estimate $\mu(x,d)\le\rho^d\mu(a,d)$ for
$x\in U$.
The claimed estimate now follows by Gronwall's lemma.
\end{pf}

Fix $R\in(0,\infty)$ and define
\[
\tilde T=\inf\{t\ge0\dvtx A_d(t)\ge RN\tilde a_{d+1}\}\wedge T.
\]

%
\begin{lemma}\label{ADDT}
There is a constant $C({\lambda},n)<\infty$ such that, for all $t\ge
0$ and all~$N$, we have
\[
\E\bigl(A_{d+1}(\tilde T\wedge t)\bigr)\le C(1+t)R^n\tilde
a_{d+1}^n+CtR^{n+1}N\tilde a_{d+1}^{n+1}.
\]
\end{lemma}
\begin{pf}
We argue as in the preceding proof, except now taking $f(x,y)=1_{\{y\ge
d+1\}}$, for which $\bar f(x)=0$.
We obtain
\[
\E\int_0^{\tilde T\wedge t}1_{\{Y_s\ge d+1\}}\,ds
\le CN^{-1}+C\E\int_0^{\tilde T\wedge t}
(Z_s^{d+1}+CN^{-1})\,ds
\]
and hence,
\begin{eqnarray*}
\E\bigl(A_{d+1}(\tilde T\wedge t)\bigr)&=&N{\lambda}\E\int_0^{\tilde T\wedge
t}(Z^{d+1}_s)^n1_{\{Y_s\ge d+1\}}\,ds\\
&\le& N{\lambda}(R\tilde a_{d+1})^n\E\int_0^{\tilde T\wedge t}1_{\{
Y_s\ge d+1\}}\,ds\\
&\le& C(1+t)R^n\tilde a_{d+1}^n+CtR^{n+1}N\tilde a_{d+1}^{n+1}.
\end{eqnarray*}
\upqed\end{pf}

\subsection{\texorpdfstring{Proof of Theorem \protect\ref{MR}}{Proof of Theorem 2.3}}
Recall that $\a$ is defined by (\ref{ALP}) and that $\a\in
(2n,2n+1)$. Note that
$\a^2-(2n+1)\a+n=0$. Recall that $\k=(2\a)^{-1}$ and
\[
d=d(N)=\sup\{k\in\N\dvtx Na_k>N^{\kappa}\}.
\]
The asymptotic growth rate (\ref{DLLN}) follows from Theorem~\ref{SGD}.
We shall use without further comment below the inequalities
\[
C^{-1}a_k^{1/\a}\le a_{k+1}\le Ca_k^\a,\qquad k\ge0,
\]
proved in Proposition~\ref{AK} and the inequalities
\[
a_{d+1}\le N^{-(1-\k)}\le a_d,\qquad d\le\log\log N,
\]
the last being valid for all sufficiently large $N$.

By the truncation estimate, Lemma~\ref{TLD}, we have
\[
\sup_{t\le t_0}\sup_{k\le d}\frac{|z_k(t)-x_k(t)|}{\sqrt{a_k}}\le
t_0\frac{a_{d+1}}{\sqrt{a_d}}\le Ct_0a_{d+1}^{1-1/(2\a)}\le Ct_0N^{-1/2}.
\]
Since $\phi(N)\to\infty$ as $N\to\infty$ it will therefore suffice
to show (\ref{MRA}) with $(z(t)\dvtx t\ge0)$ replaced by
$(x(t)\dvtx t\ge0)$.

We apply the general procedure of Section~\ref{SOTE}.
Take as norm scales $\s_k=\sqrt{a_k}$ so that
\[
\|x\|=\max_k|x_k|/\sqrt{a_k},\qquad x\in\R^d.
\]
We now identify suitable regularity constants $\L,B,\t,J,J_1(b),J(\mu),K$.
We write $C$ for a finite positive constant which may depend on
${\lambda}$ and $n$ and whose value may vary from line to line.
We shall see that, as $N\to\infty$, the inequalities between these
regularity constants required in Theorem~\ref{FLET} become valid.
The maximum jump rate is bounded above by
\[
\L=N(1+{\lambda})=CN.
\]
We refer to the form of $b(x,y)$ given at (\ref{BKX}) and note that,
for $x\in U$ and $y\in I$,
\[
\|b(x,y)\|\le B=2^na_d^{-1/2+n/\a}=Ca_d^{-1/2+n/\a}.
\]
We showed in Section~\ref{CLC} the following upper bound on the
mean coupling time of our
coupling mechanism:
\[
m(x,y,\bar y)\le\t=\frac4{N{\lambda}(1-{\lambda})^2}=CN^{-1}.
\]
We refer to Section~\ref{SOTE} for the definitions of the jump
bounds $J,J_1(b),J(\mu)$
and leave the reader to check the validity of the following inequalities:
\[
J\le N^{-1}a_d^{-1/2},\qquad J_1(b)\le CN^{-1}a_d^{-1/2+(n-1)/\a},\qquad
J(\mu)\le2nN^{-1}.
\]
Recall from (\ref{BARBE}) the form of $\bar b$.
In estimating the Lipschitz constant $K$ for $\bar b$ on $U$, first
note that, for $x\in U$ and for $j=1,\ldots,k-1$,
\[
\biggl|\frac{\pd}{\pd x_j}x_{k-1}^n\mu(x,k-1)\biggr|\le
Cx_{k-1}^n\mu(x,k-1)(x_j^{-1}+1).
\]
Here we have used the explicit form (\ref{MZK}) of $\mu(x,k-1)$ and
the fact that $(1-p_{j-1}(x))^{-1}\le C$ on $U$.
Also note the inequalities
\[
x_{k-1}^{2n-1}\sqrt{\frac{a_{k-1}}{a_k}}\le
2^{2n-1}a_{k-1}^{2n-1/2-\a/2}\le C, \qquad\sum_{j=1}^\infty\sqrt
{a_j}\le C.
\]
We find, after some further straightforward estimation, that we can
take $K=C$.

Recall the choice of function $\phi$ in the statement of Theorem \ref
{MR}. Set
\[
\ve=\sqrt{\frac{\phi(N)}N},\qquad {\delta}=\ve e^{-Kt_0}/7,\qquad {\delta
}(\b,b)={\delta},\qquad {\delta}(\g,g)={\delta}/(2\t B).
\]
Recall that $\bX_0=(1/N,0,\ldots,0)$ and $x_0=0$ and that the driving
rate $\nu$ for the coupling mechanism is equal to $N{\lambda}$.
It is now straightforward to check that
all the inequalities required in the statement of Theorem~\ref{FLE}
are valid, for all sufficiently large $N$.

Now we check the tube condition of Theorem~\ref{FLE}. The inequalities
$0\le x_d(\xi)\le\cdots\le x_1(\xi)\le1$ hold for all $\xi\in S$.
By a monotonicity property established in the proof of Theorem \ref
{DDD}, we have $x_k(t)\le a_k$ for all $t\ge0$ and for $k=1,\ldots,d$.
Hence, for $N$ sufficiently large, if $\|{\mathbf x}(\xi)-x(t)\|\le
2\ve$ for some $t\ge0$, then $x^k(\xi)\le a_k+2\ve\sqrt{a_k}\le2a_k$
and $x^1(\xi)\le a_1+2\ve\sqrt{a_1}\le{\lambda}+(1-{\lambda
})/2\le(1+{\lambda})/2$, so ${\mathbf x}(\xi)\in U$ and the tube
condition is satisfied.

Now we turn to the extra conditions needed to apply Theorem \ref
{FLET}. We noted in Section~\ref{QVB}
the quadratic variation bounds
\[
a_k(y)\le A\s_k^2,\qquad \bar a_k(x)\le\bar A\s_k^2,
\]
valid for all $x\in U$ and $y\in I$, where
\[
A=C/(Na_d^{1-n/\a}),\qquad \bar A=C(\log N)^C/N
\]
and where $\bar A\le A\le\L\bar J^2$ for sufficiently large $N$.
It is now straightforward to check, also for $N$ sufficiently large,
that the remaining
inequalities required in the statement of Theorem~\ref{FLET} hold.
Theorem~\ref{FLET} therefore applies to give
%
%
\begin{eqnarray}\label{PEST}
\PP\Bigl(\sup_{t\le t_0}\|\bX_t-x_t\|>\ve\Bigr)&\le&
2de^{-{\delta}^2/(4\bar{A}t_0)}+2de^{-(\bar{A}/A)^2t_0/(6400\Lambda
\tau^2)}\nonumber\\[-8pt]\\[-8pt]
&&{}+ \PP\bigl(\O(\beta,b)^c\cup\O(\gamma,g)^c\bigr).\nonumber
\end{eqnarray}
Now, for $N$ sufficiently large, we have $d\le\log\log N$ and, by our
choice of $\phi$ and~$\k$,
\[
{\delta}^2/(4\bar At_0)\ge\phi(N)/((\log N)^Ct_0)\ge\log N
\]
and
\[
(\bar{A}/A)^2t_0/(6400\Lambda\tau^2)\ge\log N.
\]
Hence, the first and second terms on the right-hand side of (\ref
{PEST}) tend to $0$ as \mbox{$N\to\infty$}.

Recall from (\ref{OBB}) and (\ref{OGG}) the form of the events $\O
(\b,b)$ and $\O(\g,g)$. In the present example,
the complementary exceptional events arise either as a result of
truncation or because of finite $N$ effects in the
fast variable dynamics, as shown by (\ref{BETAB}) and
estimate (\ref{GAMMAG}).
Recall that ${\delta}(\b,b)={\delta}$ and ${\delta}(\g,g)={\delta
}/(2\t B)$.
Then
%
%
\begin{equation}\label{OBBE}\quad
\O(\b,b)^c\sse\biggl\{\int_0^{T\wedge t_0}\frac{Z_t^{d+1}}{\sqrt
{a_d}}\,dt\ge{\delta}(\b,b)\biggr\}
\sse\biggl\{A_d(T\wedge t_0)\ge\frac{N{\delta}\sqrt
{a_d}}{t_0}\biggr\}.
\end{equation}
It is straightforward to check that, for all sufficiently large $N$,
${\delta}(\g,g)\ge2t_0(1+{\lambda}n)$, which implies that
%
%
\begin{eqnarray}\label{OGGE}
\O(\g,g)^c&\sse&\biggl\{\int_0^{T\wedge t_0}(1+{\lambda}n+N{\lambda
}nZ_t^{d+1})\,dt\ge{\delta}(\g,g)\biggr\}
\nonumber\\[-8pt]\\[-8pt]
&\sse&\biggl\{A_d(T\wedge t_0)\ge\frac{{\delta}}{4{\lambda}nt_0\t
B}\biggr\}.\nonumber
\end{eqnarray}
To see that $\PP(\O(\beta,b)^c\cup\O(\gamma,g)^c)\to0$ as $N\to
\infty$, we use the bound on $\E(A_d(T\wedge t_0))$
proved in Lemma~\ref{ADT} and Markov's inequality. It then suffices to
show that
in the limit $N\to\infty$,
\[
C(a_d^n+\rho^dNa_{d+1})e^{Ct_0}\ll\frac{\sqrt{\phi(N)Na_d}}{4nt_0}.
\]
For the term involving $a_d^n$ this is easy.
For the other term, involving $Na_{d+1}$, we can check that, in fact,
\[
Na_{d+1}\ll\sqrt{Na_d},\qquad \rho^d\le(\log N)^C\ll\sqrt{\phi(N)}.
\]
This completes the proof of (\ref{MRA}).
Limit (\ref{MRB}) follows immediately from Lem\-ma~\ref{ADT} using
Markov's inequality.
Finally, note that, as $N\to\infty$,
\[
\tilde a_{d+1}\le CN^{-1}a_{d}^{n}+(\log N)^Ca_{d+1}\le C\bigl(N^{-1}+(\log
N)^CN^{-1+\k}\bigr)\to0
\]
and
\[
N\tilde a_{d+1}^{n+1}\le C\bigl(N^{-(1-\k)(n+1)/\a}+(\log N)^CN^{1-(1-\k
)(n+1)}\bigr)\to0.
\]
Then the limit (\ref{MRC}) follows from (\ref{MRB}) and Lemma \ref
{ADDT} using Markov's inequality.

\subsection{Monotonicity of the queueing model}
Here we prove a natural monotonicity property of the supermarket model
with memory which is a
microscopic counterpart of the monotonicity of solutions to the
differential equation (\ref{ODE}) shown in Theorem~\ref{DDD}.
We do not rely on this result in the rest of the paper.

First we construct, on a single probability space, for all $\xi
=(z,y)\in S$,
a~version $X=X(\xi)$ of the supermarket model with memory starting
from $\xi$.
Set $y_1=y_1(\xi)=y$ and determine $y_i=y_i(\xi)\in\Z^+$ for
$i=2,\ldots,N$ by the conditions
\[
y_2\le\cdots\le y_N,\qquad z_k=\bigl|\bigl\{i\in\{1,\ldots,N\}\dvtx y_i\ge k\bigr\}\bigr|/N,
\qquad k\in\N.
\]
We work on the canonical space of a marked Poisson process of rate
$N(1+{\lambda})$, where the
marks are either, with probability $1/(1+{\lambda})$, independent
copies of a uniform random variable $J$ in
$\{1,\ldots,N\}$ or, with probability ${\lambda}/(1+{\lambda})$,
independent copies of a uniform random sample $(J_1,\ldots,J_n)$
from $\{1,\ldots,N\}$.
Fix $\xi=(z,y)\in S$ and define a process $X=X(\xi)=(X_t\dvtx t\ge0)$ in
$S$ as follows.
Set $X_t=\xi$ for all $t<T$, where $T$ is the first jump time of the
Poisson process.
If the first mark is a random variable, $J$ say, take the sequence
$y_1,\ldots,y_N$ and replace $y_J$ by $(y_J-1)^+$ to obtain a sequence
$u_1,\ldots,u_N$ say;
set $\tilde y_1=u_1$ and write $u_2,\ldots,u_N$ in nondecreasing order
to obtain $\tilde y_2\le\cdots\le\tilde y_N$.
If the first mark is a random sample, $(J_1,\ldots,J_n)$ say, select
components $(y_i\dvtx i\in\{1,J_1,\ldots,J_n\})$
and write these in nondecreasing order, $w_1\le\cdots\le w_m$ say;
replace $w_1$ by $w_1+1$ and write the resulting sequence, again in
nondecreasing order, $v_1\le\cdots\le v_m$ say;
set $\tilde y_1=v_1$ and write $v_2,\ldots,v_m$ combined with the
unselected components $(y_i\dvtx i\notin\{1,J_1,\ldots,J_n\})$
in nondecreasing order to obtain $\tilde y_2\le\cdots\le\tilde y_N$.
Set $X_T=((Z^k_T\dvtx k\in\N),Y_T)$, where
\[
Z^k_T=\bigl|\bigl\{i\in\{1,\ldots,N\}\dvtx\tilde y_i\ge k\bigr\}\bigr|/N,\qquad
k\in\N,\qquad
Y_T=\tilde y_1,
\]
and repeat the construction from $X_T$ in the usual way.

For $\xi,\xi'\in S$ write $\xi\le\xi'$ if $y_i(\xi)\le y_i(\xi
')$ for $i=1,\ldots,N$.
%
%
\begin{theorem}
Let $\xi,\xi'\in S$ with $\xi\le\xi'$. Then $X_t(\xi)\le X_t(\xi
')$ for all $t\ge0$.
\end{theorem}
\begin{pf}
It will suffice to check that the desired inequality holds at the first
jump time $T$, that is to say, with
obvious notation, that $\tilde y_i\le\tilde y'_i$ for all $i$.
Note that if $a_i\le b_i$ for all $i$ for two sequences $(a_1,\ldots
,a_n)$ and $(b_1,\ldots,b_n)$,
then the same is true for their nondecreasing rearrangements.
In the case where the first mark is a random variable $J$, since
$y_i\le y_i'$ for all $i$,
we have $u_i\le u'_i$ for all $i$
and so $\tilde y_i\le\tilde y'_i$ for all $i$.
On the other hand, when the first mark is a random sample $(J_1,\ldots
,J_n)$, we have
$w_j\le w'_j$ for all $j$, so $v_j\le v'_j$ for all $j$, and so $\tilde
y_i\le\tilde y'_i$ for all $i$.
\end{pf}


%

\printaddresses

\end{document}